\documentstyle[12pt]{article}
\def\Bbb R{{\rm \bf R}}
\def\proclaim#1{\vskip2mm{\bf #1}\em}
\def\endproclaim{\em \vskip2mm}
\def\tag#1{\eqno(#1)}
\def\gathered{\begin{array}{c}}
\def\endgathered{\end{array}}
\def\text{\mbox}

\begin{document}

\title {Integrating the probe and singular sources methods:II.  The Stokes system}
\author{Masaru IKEHATA\footnote{
Laboratory of Mathematics,
Graduate School of Advanced Science and Engineering,
Hiroshima University,
Higashihiroshima 739-8527, JAPAN.
e-mail address: ikehataprobe@gmail.com}
\footnote{Emeritus Professor at Gunma University, Maebashi 371-8510, JAPAN}
}
\maketitle

\begin{abstract}
In this paper, an integrated theory of the probe and singular sources methods
for an inverse obstacle problem governed by the Stokes system in a bounded domain is developed.
The main results consist of:
the probe method for the Stokes system;
the singular sources method by using the notion of the probe method;
the completely integrated version of the probe and singular sources methods.
In establishing the singular sources method, a third indicator function which is called
the IPS function plays the central role.

\noindent
AMS: 35E05, 35J05, 35R30, 76D07 

\end{abstract}


\section{Introduction}

As classical analytical methods for inverse obstacle problems governed by partial differential equations, the probe method of Ikehata \cite{IProbe} and singular sources method of Potthast \cite{P1} are well-known.
Recently in \cite{IPS}  the author introduced a theory which integrates
the probe and singular sources methods by considering a prototype inverse obstacle problem
governed by the Laplace equation.

The aim of this paper is to extend the range of this theory which we call the {\it integrated theory} of 
the probe and singular sources methods (IPS)
to an inverse obstacle problem governed by the Stokes system in a bounded domain.

\subsection{Skelton of IPS}

First let us explain about the framework of the {\it integrated theory} of the probe and singular sources methods
developed in \cite{IPS}.
To do this, it would be a good idea to start with a description of the prototypical inverse obstacle problem considered therein.

Let $\Omega$ be a bounded domain of $\Bbb R^3$ with smooth boundary.
Let $D$ be a non empty open set of $\Bbb R^3$ with smooth boundary satisfying that
$\overline{D}\subset\Omega$
and $\Omega\setminus\overline{D}$ is connected.   The set $D$ is a mathematical model of  a cavity occurred
inside the background domain $\Omega$.

Given a function $f$ on $\partial\Omega$ let $u$ be the solution
of the boundary value problem
$$
\left\{
\begin{array}{ll}
\displaystyle
\Delta u(z)=0, & z\in\Omega\setminus\overline{D},
\\
\\
\displaystyle
\frac{\partial u}{\partial\nu}(z)=0, & z\in\partial D,
\\
\\
\displaystyle
u(z)=f(z), & z\in\partial\Omega.
\end{array}
\right.
$$
Here we denote by $\nu$ the unit outer normal on $\partial\Omega$ and $\partial D$.

The inverse obstacle problem considered in \cite{IPS} is:
reconstruct $D$ from the Dirichlet-to-Neumann map
$$\displaystyle
\Lambda_D:f\longmapsto\frac{\partial u}{\partial\nu}\vert_{\partial\Omega}.
$$
We denote $\Lambda_D$ by $\Lambda_0$ if $D=\emptyset$.

Let $x\in\Omega$ be an arbitrary point in $\Omega$ and 
denote by $G(z,x)$ the standard fundamental solution of the Laplace equation, that is
$$\displaystyle
G(z,x)=\frac{1}{4\pi\vert z-x\vert}.
$$
The center of IPS is the function $W(z)=W_x(z)$ of  independent variables $z\in\Omega\setminus\overline{D}$
for each fixed $x\in\Omega\setminus\overline{D}$ as the solution 
of the boundary value problem
$$
\left\{\begin{array}{ll}
\Delta W(z)=0, & z\in\Omega\setminus\overline{D},
\\
\\
\displaystyle
\frac{\partial W}{\partial\nu}(z)=-\frac{\partial}{\partial\nu}G(z,x), & z\in\partial D,
\\
\\
\displaystyle
W(z)=G(z,x),  & z\in\partial\Omega.
\end{array}
\right.
\tag {1.1}
$$
Then the {\it main} discovery of \cite{IPS} is concerned with the function $W_x(x)\equiv W_x(z)\vert_{z=x}$ for $x\in\Omega\setminus\overline{D}$, which is called the {\it third indicator function}
or the IPS function.
Since the governing equation is linear, one has the natural decomposition
$$\displaystyle
W_x(x)=w_x(x)+w_x^1(x),
\tag {1.2}
$$
where
$w_x(x)\equiv w_x(z)\vert_{z=x}$, $w_x^1(x)=w_x^1(z)\vert_{z=x}$ and 
the functions $w_x$ and $w_x^1$ satisfy
$$
\left\{\begin{array}{ll}
\Delta w_x(z)=0, & z\in\Omega\setminus\overline{D},
\\
\\
\displaystyle
\frac{\partial w_x}{\partial\nu}(z)=-\frac{\partial}{\partial\nu}G(z,x), & z\in\partial D,
\\
\\
\displaystyle
w_x(z)=0,  & z\in\partial\Omega
\end{array}
\right.
\tag {1.3}
$$
and
$$
\left\{\begin{array}{ll}
\Delta w_x^1(z)=0, & z\in\Omega\setminus\overline{D},
\\
\\
\displaystyle
\frac{\partial w_x^1}{\partial\nu}(z)=0, & z\in\partial D,
\\
\\
\displaystyle
w_x^1(z)=G(z,x),  & z\in\partial\Omega.
\end{array}
\right.
\tag {1.4}
$$
In \cite{IPS} we have found another decomposition of IPS function, which is called the {\it inner decomposition}, that is
$$\displaystyle
W_x(x)=I(x)+I^1(x),
\tag {1.5}
$$
where the functions $I(x)$ and $I^1(x)$ are given by
$I(x)\equiv I(x,y)\vert_{x=y}$ and $I^1(x,y)\vert_{x=y}$ and
$$
\left\{
\begin{array}{l}
\displaystyle
I(x,y)=\int_{\Omega\setminus\overline{D}}\nabla w_x(z)\cdot\nabla w_y(z)\,dz
+\int_D\nabla G(z,x)\cdot\nabla G(z,y)\,dz,
\\
\\
\displaystyle
I^1(x,y)=
\int_{\Omega\setminus\overline{D}}\nabla w_x^1(z)\cdot\nabla w_y^1(z)\,dz
+\int_{\Bbb R^3\setminus\overline{D}}
\nabla G(z,x)\cdot\nabla G(z,y)\,dz.
\end{array}
\right.
$$

We see that  the $I(x)$ of $x\in\Omega\setminus\overline{D}$ in the inner decomposition (1.5) coincides with the indicator function of the probe method in \cite{INew}.
More precisely, recall the notion of needle reformulated in \cite{INew}.
Given $x\in\Omega$ a non self-intersecting piecewise linear curve $\sigma=\sigma(t)$, $0\le t\le 1$ is called a needle
with a tip at $x$ if $\sigma(0)\in\partial\Omega$, $\sigma(1)=x$ and $\sigma(t)\in\Omega$ for all $t\in\,]0,\,1[$.
We denote by $N_x$ the set of all needles with a tip at $x$.
Then we have the expression
$$\displaystyle
I(x)=\lim_{n\rightarrow\infty}
\int_{\partial\Omega}(\Lambda_0-\Lambda_D)(v_n\vert_{\partial\Omega})(z)\,
v_n(z)\,dS(z),
$$
where $\{v_n\}$ is a sequence of harmonic functions in $\Omega$ that converges to $G(\,\cdot\,,x)$ in $\Omega\setminus\sigma$ compact uniformly for a needle $\sigma\in N_x$
and satisfying $\sigma\cap\overline{D}=\emptyset$.
And the blowing  up property
$$\displaystyle
\lim_{x\rightarrow a\in\partial D} I(x)=\infty,
$$
is valid.

The points of  IPS are as follows.

$\quad$

\noindent
$\bullet$ The $w_x(x)$ of $x\in\Omega\setminus\overline{D}$ in the natural decomposition (1.2)
also has the expression
$$\displaystyle
w_x(x)=-\lim_{n\rightarrow\infty}
\int_{\partial\Omega}(\Lambda_0-\Lambda_D)(v_n\vert_{\partial\Omega})(z)\,
(G(z,x)-v_n(z))\,dS(z)
$$
and the blowing  up property
$$\displaystyle
\lim_{x\rightarrow a\in\partial D} w_x(x)=\infty.
$$
This blowing up property is deduced from that of $I(x)$ by using the decompositions (1.2) and (1.5)
without analyzing $w_x(y)$ itself.

$\quad$

\noindent
$\bullet$  The functions $W_x(y)$ together with $w_x(y)$, $w_x^1(y)$, $I(x,y)$ and $I^1(x,y)$ can be calculated from
$\Lambda_D$ and $\Lambda_0$ with or without using the needle sequences.  The IPS function
$W_x(x)$ blows up on both $\partial\Omega$ and $\partial D$.  The functions
$w_x^1(x)$ and $I_x^1(x)$ are bounded in a neighbourhood of $\partial D$.

$\quad$

\noindent
$\bullet$  In \cite{INew} we have already known that the sequence, which is called the indicator sequence for the probe method
given by
$$\displaystyle
\left\{\int_{\partial\Omega}(\Lambda_0-\Lambda_D)(v_n\vert_{\partial\Omega})(z)\,
v_n(z)\,dS(z)\right\},
$$
yields the reconstruction formula of $\overline{D}$ itself by its blowing up property.
In \cite{IPS} in addition to this, we have found another sequence playing the same role as the indicator sequence, that is,
the sequence
$$\displaystyle
\left\{\int_{\partial\Omega}(\Lambda_0-\Lambda_D)((G(\,\cdot\,,x)-v_n)\vert_{\partial\Omega})(z)\,
(G(z,x)-v_n(z))\,dS(z)\right\},
$$
gives  us the reconstruction formula of $\overline{D}$ itself by its blowing up property tested
for all $x\in\Omega$, $\sigma\in N_x$ and 
sequences of harmonic functions in $\Omega$ that converge to $G(\,\cdot\,,x)$ in $\Omega\setminus\sigma$ compact uniformly.
So we have two reconstruction formulae of $\overline{D}$.  

$\quad$

\noindent
$\bullet$  And as a byproduct, we found the completely integrated version of the two methods, in the sense that
a single indicator function together with its corresponding indicator sequence plays the role, or meaning of those of the probe and singular sources methods
at the same time.
The idea is just replace $G(z,x)$ in (1.1), (1.3), (1.4) with the Green function $G_{\Omega}(z,x)$ for the Laplace equation in the whole domain $\Omega$.

$\quad$

\noindent
This is the skelton of the integrated theory of the probe and singular sources methods.

Note that the original probe method has been applied to various inverse obstacle problems governed by, for example,
the conductivity equation coming from the Calder\'on problem \cite{Cal}, the Navier equation from the nondestructive evaluation
of the material, Helmholtz equation $\Delta u+k^2u=0$ coming from inverse `obstacle' scattering problems in a wider sense.  See Section 2 of the author's review paper \cite{IReview} for past works and Section 5.2 of \cite{IRevisit}
on recent progress.
The original singular sources method also has been widely applied, starting from 
inverse obstacle scattering problems governed by the Helmholtz equation.  See the book \cite{NP} and references therein.

However, as far as the author knows, as the founder of the probe method, 
both methods has not yet been theoretically realized 
in inverse obstacle problems governed by the Stokes system in the bounded domain compared with applications
\cite{LR}, \cite{LY}, \cite{MS}, \cite{HUW}
of other classical methods belonging to the same group, such as
the factorization method of Kirsh \cite{K}, the linear sampling method of Colton-Kirsh \cite{CKi}
and the enclosure method of Ikehata \cite{IEn}.
Besides, it is natural to ask whether the integrated theory of the probe and singular sources 
methods can be realized in other interesting inverse problems governed by various partial differential equations
in the bounded domain.

In this paper we choose a problem governed by the Stokes system.

\subsection{Problem formulation}

Let $\Omega\subset\Bbb R^3$ be an arbitrary bounded domain of $\Bbb R^3$ with $C^2$-boundary.
Let $D$ be a nonempty open set of $\Bbb R^3$ such that $\overline{D}\subset\Omega$,
$\Omega\setminus\overline{D}$ is connected and $\partial D$ is $C^2$.
In what follows the $\mbox{\boldmath $\nu$}$ always denotes the unit outward normal vector field on $\partial\Omega$ and $\partial D$.

Given $\mbox{\boldmath $f$}\in H^{\frac{3}{2}}(\partial\Omega)$ satisfying $\int_{\partial\Omega}\mbox{\boldmath $f$}\cdot\mbox{\boldmath $\nu$}\,dS=0$, let
 $(\mbox{\boldmath $u$},p)\in H^2(\Omega\setminus\overline{D})\times H^1(\Omega\setminus\overline{D})$
 the solution of  the Dirichlet problem for the Stokes system
 $$\left\{
 \begin{array}{ll}
 \mu\Delta\mbox{\boldmath $u$}-\nabla p=\mbox{\boldmath $0$}, & x\in\Omega\setminus\overline{D},
 \\
 \\
 \displaystyle
 \nabla\cdot\mbox{\boldmath $u$}=0, &  x\in\Omega\setminus\overline{D},\\
 \\
 \displaystyle
 \mbox{\boldmath $u$}=\mbox{\boldmath $0$}, & x\in\partial D,
 \\
 \\
 \displaystyle
 \mbox{\boldmath $u$}=\mbox{\boldmath $f$}, & x\in\partial\Omega.
 \end{array}
 \right.
 \tag {1.6}
 $$
Here $\mu$ is a positive constant and denotes the viscosity of the fluid occupying $\Omega\setminus\overline{D}$,
and the pair $(\mbox{\boldmath $u$},p)$  denotes that of the velocity and pressure of the fluid.

Note that, in this paper, standard results regarding the well-posedness
(the existence, uniqueness, stability) of the Stokes system in the Sobolev spaces are used.
See, for example, the facts  summarized in Theorem 2.2 of \cite{BA}.  In particular, both the $\mbox{\boldmath $u$}$ and $p$ exist;
$\mbox{\boldmath $u$}$ is unique and so is $p$, however, modulo a constant.
More precisely, we use only the following facts about the direct problem extracted from Theorem 2.2 in \cite{BA}.

Given $(\mbox{\boldmath $g$}_1,\mbox{\boldmath $g$}_2)\in H^{\frac{3}{2}}(\partial\Omega)\times
H^{\frac{3}{2}}(\partial D)$ with $\int_{\partial\Omega}\mbox{\boldmath $g$}_1\cdot\mbox{\boldmath $\nu$}\,dS=0$
and $\int_{\partial D}\mbox{\boldmath $g$}_2\cdot\mbox{\boldmath $\nu$}\,dS=0$,
there exists a solution $(\mbox{\boldmath $u$},p)\in H^{2}(\Omega\setminus\overline{D})
\times H^1(\Omega\setminus\overline{D})$ such that
$$\left\{
 \begin{array}{ll}
 \mu\Delta\mbox{\boldmath $u$}-\nabla p=\mbox{\boldmath $F$}, & x\in\Omega\setminus\overline{D},
 \\
 \\
 \displaystyle
 \nabla\cdot\mbox{\boldmath $u$}=0, &  x\in\Omega\setminus\overline{D},\\
 \\
 \displaystyle
 \mbox{\boldmath $u$}=\mbox{\boldmath $g$}_2, & x\in\partial D,
 \\
 \\
 \displaystyle
 \mbox{\boldmath $u$}=\mbox{\boldmath $g$}_1, & x\in\partial\Omega
 \end{array}
 \right.
 $$
 and
 $$\displaystyle
 \Vert\mbox{\boldmath $u$}\Vert_{H^2(\Omega\setminus\overline{D})}
 +\Vert p-p_{E}\Vert_{H^1(\Omega\setminus\overline{D})}
 \le C\,(\Vert\mbox{\boldmath $g$}_1\Vert_{H^{\frac{3}{2}}(\partial\Omega)}+
 \Vert\mbox{\boldmath $g$}_1\Vert_{H^{\frac{3}{2}}(\partial D)}+\Vert\mbox{\boldmath $F$}\Vert_{L^2(\Omega\setminus\overline{D})}),
 $$
 where $C$ is a positive constant independent of the Dirichlet data $\mbox{\boldmath $g$}_1$ and $\mbox{\boldmath $g$}_2$ and $P_{E}$ with $E=\Omega\setminus\overline{D}$ denotes
 $$\displaystyle
 p_{E}=\frac{\int_{E} p\,dx}{\vert E\vert}.
 $$
 In this paper, we do not cite these facts one by one, but use the term `well-posedness' for the Stokes system.

 $\quad$

{\bf\noindent Definition 1.1.}  Given $\mbox{\boldmath $f$}\in H^{\frac{3}{2}}(\partial\Omega)$ with
$\int_{\partial\Omega}\mbox{\boldmath $f$}\cdot\mbox{\boldmath $\nu$}\,dS=0$, let
$(\mbox{\boldmath $u$},p)$ be the solution of (1.6).
Define the Dirichlet-to-Neumann map $\Lambda_D$ by the formula
$$\begin{array}{ll}
\displaystyle
\Lambda_D\mbox{\boldmath $f$}(x)=
\sigma(\mbox{\boldmath $u$},p)\mbox{\boldmath $\nu$}(x), & x\in\partial\Omega.
\end{array}
$$
where
$$\displaystyle
\sigma(\mbox{\boldmath $u$},p)
=2\mu\,\text{Sym}\nabla\mbox{\boldmath $u$}(x)-p(x)\,I_3.
$$

$\quad$

\noindent
Note that our observation date is the stress $\sigma(\mbox{\boldmath $u$},p)\mbox{\boldmath $\nu$}$ 
distributed on $\partial\Omega$ itself that induces the velocity field $\mbox{\boldmath $f$}$ on $\partial\Omega$.

The relationship of the stress tensor $\sigma(\mbox{\boldmath $u$},p)$ and the original equation is as follows:
$$\begin{array}{ll}
\displaystyle
\{\text{div}\,\sigma(\mbox{\boldmath $u$},p)\}\cdot\mbox{\boldmath $a$}
&
\displaystyle
=
\nabla\cdot(\sigma(\mbox{\boldmath $u$},p)\mbox{\boldmath $a$})\\
\\
\displaystyle
&
\displaystyle
=\{\mu\Delta\mbox{\boldmath $u$}+\mu\nabla(\nabla\cdot\mbox{\boldmath $u$})-\nabla p\}\cdot\mbox{\boldmath $a$},
\end{array}
$$
where $\mbox{\boldmath $a$}$ is an arbitrary constant vector.
Thus under the condition $\nabla\cdot\mbox{\boldmath $u$}=0$, the equation
$$\displaystyle
\text{div}\,\sigma(\mbox{\boldmath $u$},p)=\mbox{\boldmath $0$}
$$
coincides with the equation 
$$\displaystyle
\mu\Delta\mbox{\boldmath $u$}-\nabla p=\mbox{\boldmath $0$}.
$$

Let $\mbox{\boldmath $\phi$}\in H^{1}(\Omega\setminus\overline{D})$.
Integration by parts yields
$$\begin{array}{l}
\displaystyle
\,\,\,\,\,\,
\int_{\partial\Omega}\Lambda_D\mbox{\boldmath $f$}\cdot\mbox{\boldmath $\phi$}\,dS
\\
\\
\displaystyle
=\int_{\Omega\setminus\overline{D}}
\text{div}\,\sigma(\mbox{\boldmath $u$},p)\cdot\mbox{\boldmath $\phi$}\,dx
+\int_{\partial D}\sigma(\mbox{\boldmath $u$},p)\mbox{\boldmath $\nu$}\cdot\mbox{\boldmath $\phi$}\,dS
+\int_{\Omega\setminus\overline{D}}\sigma(\mbox{\boldmath $u$},p)\cdot\nabla\mbox{\boldmath $\phi$}\,dx
\\
\\
\displaystyle
=\int_{\partial D}\sigma(\mbox{\boldmath $u$},p)\mbox{\boldmath $\nu$}\cdot\mbox{\boldmath $\phi$}\,dS+\int_{\Omega\setminus\overline{D}}\sigma(\mbox{\boldmath $u$},p)\cdot\nabla\mbox{\boldmath $\phi$}\,dx\\
\\
\displaystyle
=\int_{\partial D}\sigma(\mbox{\boldmath $u$},p)\mbox{\boldmath $\nu$}\cdot\mbox{\boldmath $\phi$}\,dS+\int_{\Omega\setminus\overline{D}}2\mu\text{Sym}\,\nabla\mbox{\boldmath $u$}\cdot
\text{Sym}\,\nabla\mbox{\boldmath $\phi$}\,dx
-\int_{\Omega\setminus\overline{D}}p(x)\nabla\cdot\mbox{\boldmath $\phi$}\,dx.
\end{array}
$$
So if $\nabla\cdot\mbox{\boldmath $\phi$}=0$ in $\Omega\setminus\overline{D}$, then we have the expression
$$\displaystyle
\int_{\partial\Omega}\Lambda_D\mbox{\boldmath $f$}\cdot\mbox{\boldmath $\phi$}\,dS
=\int_{\partial D}\sigma(\mbox{\boldmath $u$},p)\mbox{\boldmath $\nu$}\cdot\mbox{\boldmath $\phi$}\,dS+\int_{\Omega\setminus\overline{D}}2\mu\,\text{Sym}\,\nabla\mbox{\boldmath $u$}\cdot
\text{Sym}\,\nabla\mbox{\boldmath $\phi$}\,dx.
\tag {1.7}
$$
Let the pair $(\mbox{\boldmath $v$},q)\in H^2(\Omega)\times H^1(\Omega)$ solve the Stokes system 
$$\left\{
 \begin{array}{ll}
 \mu\Delta\mbox{\boldmath $v$}-\nabla q=\mbox{\boldmath $0$}, & x\in\Omega,
 \\
 \\
 \displaystyle
 \nabla\cdot\mbox{\boldmath $v$}=0, & x\in\Omega,
 \\
 \\
 \displaystyle
 \mbox{\boldmath $v$}=\mbox{\boldmath $f$}, & x\in\partial\Omega.
 \end{array}
 \right.
 $$
We  denote by $\Lambda_D=\Lambda_0$ in the case when $D=\emptyset$.
In this case we have the expression
$$\displaystyle
\int_{\partial\Omega}\Lambda_0\mbox{\boldmath $f$}\cdot\mbox{\boldmath $\Psi$}\,dS
=\int_{\Omega}2\mu\,\text{Sym}\,\nabla\mbox{\boldmath $v$}\cdot
\text{Sym}\,\nabla\mbox{\boldmath $\Psi$}\,dx,
\tag {1.8}
$$
where $\mbox{\boldmath $\Psi$}\in H^1(\Omega)$ with $\nabla\cdot\mbox{\boldmath $\Psi$}=0$ in $\Omega$.

Here we derive an important representation formula of $\Lambda_0-\Lambda_D$.

\proclaim{\noindent Proposition 1.1.}
Let $(\mbox{\boldmath $v$},q)=(\mbox{\boldmath $v$}_j,q_j)\in H^2(\Omega)\times H^1(\Omega)$, $j=1,2$ be arbitrary  solutions of
$$\left\{
 \begin{array}{ll}
 \mu\Delta\mbox{\boldmath $v$}-\nabla q=\mbox{\boldmath $0$}, & x\in\Omega,
 \\
 \\
 \displaystyle
 \nabla\cdot\mbox{\boldmath $v$}=0, & x\in\Omega.
 \end{array}
 \right.
 \tag {1.9}
 $$
And  let $(\mbox{\boldmath $u$},p)=(\mbox{\boldmath $u$}_j,p_j)$, $j=1,2$ solve (1.6) with $\mbox{\boldmath $f$}=\mbox{\boldmath $v$}_j$, $j=1,2$, respectively.

\noindent
Then we have
$$\begin{array}{l}
\displaystyle
\,\,\,\,\,\,
\int_{\partial\Omega}(\Lambda_0-\Lambda_D)\mbox{\boldmath $v$}_1
\vert_{\partial\Omega}\cdot\mbox{\boldmath $v$}_2\vert_{\partial\Omega}\,dS
\\
\\
\displaystyle
=-\int_{\Omega\setminus\overline{D}}2\mu\,
\text{Sym}\,\nabla\mbox{\boldmath $w$}_1\cdot
\text{Sym}\,\nabla\mbox{\boldmath $w$}_2\,dx
-\int_{D}2\mu\,
\text{Sym}\,\nabla\mbox{\boldmath $v$}_1\cdot
\text{Sym}\,\nabla\mbox{\boldmath $v$}_2\,dx,
\end{array}
\tag {1.10}
$$
where $\mbox{\boldmath $w$}_j=\mbox{\boldmath $u$}_j-\mbox{\boldmath $v$}_j$, $j=1,2$.

\endproclaim

{\it\noindent Proof.}
Noting $\mbox{\boldmath $v$}_1=\mbox{\boldmath $u$}_1$ and $\mbox{\boldmath $v$}_2=\mbox{\boldmath $u$}_2$ on $\partial\Omega$, and $\mbox{\boldmath $u$}_2=\mbox{\boldmath $0$}$ on $\partial D$, from (1.7) and (1.8) we have
$$\begin{array}{l}
\displaystyle
\,\,\,\,\,\,
\int_{\partial\Omega}(\Lambda_0-\Lambda_D)\mbox{\boldmath $v$}_1
\vert_{\partial\Omega}\cdot\mbox{\boldmath $v$}_2\vert_{\partial\Omega}\,dS
\\
\\
\displaystyle
=\int_{\partial\Omega}\Lambda_0\mbox{\boldmath $v$}_1\vert_{\partial\Omega}\cdot\mbox{\boldmath $v$}_2\vert_{\partial\Omega}dS
-\int_{\partial\Omega}\Lambda_D\mbox{\boldmath $u$}_1\vert_{\partial\Omega}\cdot\mbox{\boldmath $u$}_2\vert_{\partial\Omega}dS
\\
\\
\displaystyle
=\int_{\Omega}2\mu\,\text{Sym}\,\nabla\mbox{\boldmath $v$}_1\cdot
\text{Sym}\,\nabla\mbox{\boldmath $v$}_2\,dx
-\int_{\Omega\setminus\overline{D}}2\mu\,\text{Sym}\,\nabla\mbox{\boldmath $u$}_1\cdot
\text{Sym}\,\nabla\mbox{\boldmath $u$}_2\,dx
\\
\\
\displaystyle
=\int_{D}2\mu\,\text{Sym}\,\nabla\mbox{\boldmath $v$}_1\cdot
\text{Sym}\,\nabla\mbox{\boldmath $v$}_2\,dx
+\int_{\Omega\setminus\overline{D}}2\mu\,\text{Sym}\,\nabla\mbox{\boldmath $v$}_1\cdot
\text{Sym}\,\nabla\mbox{\boldmath $v$}_2\,dx
\\
\\
\displaystyle
\,\,\,
-\int_{\Omega\setminus\overline{D}}2\mu\,
\text{Sym}\,\nabla(\mbox{\boldmath $w$}_1+\mbox{\boldmath $v$}_1)\cdot
\text{Sym}\,\nabla(\mbox{\boldmath $w$}_2+\mbox{\boldmath $v$}_2)\,dx\\
\\
\displaystyle
=\int_{D}2\mu\,\text{Sym}\,\nabla\mbox{\boldmath $v$}_1\cdot
\text{Sym}\,\nabla\mbox{\boldmath $v$}_2\,dx
\\
\\
\displaystyle
\,\,\,
-\int_{\Omega\setminus\overline{D}}2\mu\,
\text{Sym}\,\nabla\mbox{\boldmath $w$}_1\cdot
\text{Sym}\,\nabla\mbox{\boldmath $w$}_2\,dx
\\
\\
\displaystyle
\,\,\,
-\int_{\Omega\setminus\overline{D}}2\mu\,
\text{Sym}\,\nabla\mbox{\boldmath $w$}_1\cdot
\text{Sym}\,\nabla\mbox{\boldmath $v$}_2\,dx
-\int_{\Omega\setminus\overline{D}}2\mu\,
\text{Sym}\,\nabla\mbox{\boldmath $v$}_1\cdot
\text{Sym}\,\nabla\mbox{\boldmath $w$}_2\,dx.
\end{array}
\tag {1.11}
$$
Here we have
$$\begin{array}{ll}
\displaystyle
\int_{\Omega\setminus\overline{D}}2\mu\,
\text{Sym}\,\nabla\mbox{\boldmath $w$}_1\cdot
\text{Sym}\,\nabla\mbox{\boldmath $v$}_2\,dx
&
\displaystyle
=\int_{\Omega\setminus\overline{D}}\sigma(\mbox{\boldmath $v$}_2, q_2)\cdot\nabla\mbox{\boldmath $w$}_1\,dx\\
\\
\displaystyle
&
\displaystyle
=-\int_{\partial D}\sigma(\mbox{\boldmath $v$}_2, q_2)\mbox{\boldmath $\nu$}\cdot\mbox{\boldmath $w$}_1\,dS
\\
\\
\displaystyle
&
\displaystyle
=\int_{\partial D}\sigma(\mbox{\boldmath $v$}_2, q_2)\mbox{\boldmath $\nu$}\cdot\mbox{\boldmath $v$}_1\,dS
\\
\\
\displaystyle
&
\displaystyle
=\int_D \text{div}\,\sigma(\mbox{\boldmath $v$}_2,q_2)\cdot\mbox{\boldmath $v$}_1\,dx
+\int_{D}\,\sigma(\mbox{\boldmath $v$}_2,q_2)\cdot\nabla\mbox{\boldmath $v$}_1\,dx
\\
\\
\displaystyle
&
\displaystyle
=\int_{D}2\mu\,\text{Sym}\,\nabla\mbox{\boldmath $v$}_2\cdot\text{Sym}\,\nabla\mbox{\boldmath $v$}_1\,dx
\end{array}
$$
and also
$$
\displaystyle
\int_{\Omega\setminus\overline{D}}2\mu\,
\text{Sym}\,\nabla\mbox{\boldmath $v$}_1\cdot
\text{Sym}\,\nabla\mbox{\boldmath $w$}_2\,dx
=\int_{D}2\mu\,\text{Sym}\,\nabla\mbox{\boldmath $v$}_1\cdot\text{Sym}\,\nabla\mbox{\boldmath $v$}_2\,dx.
$$
Substituting these into (1.11), we obtain the identity (1.10).

\noindent
$\Box$

\noindent
Note that the reason for the non-appearance of
the pressures  on the right-hand side of (1.10)  is: we have made use of the divergence free property:
$$\left\{
\begin{array}{ll}
\displaystyle
\nabla\cdot\mbox{\boldmath $u$}_j=0, & y\in\Omega\setminus\overline{D},\\
\\
\displaystyle
\nabla\cdot\mbox{\boldmath $v$}_j=0, & y\in\Omega.
\end{array}
\right.
$$

From (1.10) and the well-poseness for the Stokes system we have the symmetry:
for all $\mbox{\boldmath $f$}\in H^{\frac{3}{2}}(\partial\Omega)$ and $\mbox{\boldmath $g$}\in H^{\frac{3}{2}}
(\partial\Omega)$ with $\int_{\partial\Omega}\mbox{\boldmath $f$}\cdot\mbox{\boldmath $\nu$}\,dS=
\int_{\partial\Omega}\mbox{\boldmath $g$}\cdot\mbox{\boldmath $\nu$}\,dS=0$, it holds that
$$\displaystyle
\int_{\partial\Omega}\,(\Lambda_0-\Lambda_D)\mbox{\boldmath $f$}\cdot\mbox{\boldmath $g$}\,dS
=\int_{\partial\Omega}\,(\Lambda_0-\Lambda_D)\mbox{\boldmath $g$}\cdot\mbox{\boldmath $f$}\,dS.
$$

In this paper, we consider the inverse obstacle problem  formulated as

$\quad$

{\it\noindent\bf Problem 1.1.}
 Extract information about the geometry of $D$ from the $\Lambda_D$.

 $\quad$
 
 In particular, we are interested in developing the integrated theory of the probe and singular sources methods
 \cite{IPS} to the Stokes system.
 Thus the mathematical model of the observation data used in this paper is the integral
 $$\displaystyle
 \int_{\partial\Omega}(\Lambda_0-\Lambda_D)\mbox{\boldmath $v$}\vert_{\partial\Omega}
 \cdot\mbox{\boldmath $v$}\vert_{\partial\Omega}\,dS
 \equiv <(\Lambda_0-\Lambda_D)\mbox{\boldmath $v$}\vert_{\partial\Omega}, 
 \mbox{\boldmath $v$}\vert_{\partial\Omega}>,
 \tag {1.12}
 $$
 where $\mbox{\boldmath $v$}\in H^2(\Omega)$ satisfies the Stokes system (1.9)
 with a $q\in H^1(\Omega)$.

\section{Singular solution, needle and needle sequence}

\subsection{The Stokeslet as the singular solution}

Let $\mbox{\boldmath $a$}$ be  an arbitrary nonzero constant vector.
It is well known that the fields $\mbox{\boldmath $u$}$ and $r$ given by
$$
\left\{\begin{array}{l}
\displaystyle
\mbox{\boldmath $u$}(x)
=\mbox{\boldmath $J$}(x)\mbox{\boldmath $a$},\\
\\
\displaystyle
r=\mbox{\boldmath $p$}(x)\cdot\mbox{\boldmath $a$}
\end{array}
\right.
\tag {2.1}
$$
satisfy the Stokes system
in the whole space $\Bbb R^3$
$$
\left\{
\begin{array}{l}
\displaystyle
\mu\Delta\mbox{\boldmath $u$}-\nabla r+\mbox{\boldmath $a$}\,\delta(x)
=\mbox{\boldmath $0$},\\
\\
\displaystyle
\nabla\cdot\mbox{\boldmath $u$}=\mbox{\boldmath $0$},
\end{array}
\right.
\tag {2.2}
$$
where the matrix valued distribution $\mbox{\boldmath $J$}(x)$ and vector valued distribution $\mbox{\boldmath $p$}(x)$ in the whole space $\Bbb R^3$ are given by
$$\left\{
\begin{array}{l}
\displaystyle
\mbox{\boldmath $J$}(x)=
\frac{1}{8\pi\mu}
\left(\frac{1}{\vert x\vert}I_3+\frac{x\otimes x}{\vert x\vert^3}\,\right),
\\
\\
\displaystyle
\mbox{\boldmath $p$}(x)=\frac{x}{4\pi\vert x\vert^3}=-\nabla\left(\frac{1}{4\pi\vert x\vert}\right).
\end{array}
\right.
$$
Those are called the Stokeslet.

The construction of the Stokeslet is a consequence of the following fact.
\proclaim{\noindent Proposition 2.1.}
Given $\mbox{\boldmath $F$}$ let $\mbox{\boldmath $\phi$}$ be a solution of the vector Laplace equation
$$\displaystyle
\Delta\mbox{\boldmath $\phi$}-\mbox{\boldmath $F$}=\mbox{\boldmath $0$}.
$$
Define the vector and scalar fields by the formula
$$
\left\{
\begin{array}{l}
\displaystyle
\mbox{\boldmath $v$}(x)=\nabla(x\cdot\mbox{\boldmath $\phi$})-2\mbox{\boldmath $\phi$},
\\
\\
\displaystyle
q(x)=2\mu\nabla\cdot\mbox{\boldmath $\phi$}.
\end{array}
\right.
$$
Then we have
$$
\left\{
\begin{array}{l}
\displaystyle
\mu\Delta\mbox{\boldmath $v$}-\nabla q+\mu\mbox{\boldmath $F$}-\mu(\nabla\mbox{\boldmath $F$})^Tx=\mbox{\boldmath $0$},
\\
\\
\displaystyle
\nabla\cdot\mbox{\boldmath $v$}=x\cdot\mbox{\boldmath $F$}.
\end{array}
\right.
$$
\endproclaim

{\it\noindent Proof.}
We have
$$\begin{array}{ll}
\displaystyle
\nabla\cdot\mbox{\boldmath $v$}
&
\displaystyle
=\Delta(x\cdot\mbox{\boldmath $\phi$})-2\nabla\cdot\mbox{\boldmath $\phi$}
\\
\\
\displaystyle
&
\displaystyle
=2\nabla\cdot\mbox{\boldmath $\phi$}+x\cdot\Delta\mbox{\boldmath $\phi$}-2\nabla\cdot\mbox{\boldmath $\phi$}
\\
\\
\displaystyle
&
\displaystyle
=x\cdot\mbox{\boldmath $F$},
\end{array}
$$
$$\displaystyle
\nabla(\nabla\cdot\mbox{\boldmath $v$})
=\mbox{\boldmath $F$}+(\nabla\mbox{\boldmath $F$})^Tx
$$
and
$$\begin{array}{ll}
\displaystyle
\Delta\mbox{\boldmath $v$}
&
\displaystyle
=\nabla(\Delta(x\cdot\mbox{\boldmath $\phi$}))-2\Delta\mbox{\boldmath $\phi$}
\\
\\
\displaystyle
&
\displaystyle
=\nabla(2\nabla\cdot\mbox{\boldmath $\phi$})+x\cdot\Delta\mbox{\boldmath $\phi$})-2\Delta\mbox{\boldmath $\phi$}
\\
\\
\displaystyle
&
\displaystyle
=2\nabla(\nabla\cdot\mbox{\boldmath $\phi$})+\Delta\mbox{\boldmath $\phi$}+\{\nabla(\Delta\mbox{\boldmath $\phi$})\}^Tx-2\Delta\mbox{\boldmath $\phi$}
\\
\\
\displaystyle
&
\displaystyle
=2\nabla(\nabla\cdot\mbox{\boldmath $\phi$})-\mbox{\boldmath $F$}
+(\nabla\mbox{\boldmath $F$})^Tx.
\end{array}
$$
Thus
$$\displaystyle
\mu\Delta\mbox{\boldmath $v$}
-\nabla(2\mu\nabla\cdot\mbox{\boldmath $\phi$})+\mu\mbox{\boldmath $F$}-\mu(\nabla\mbox{\boldmath $F$})^Tx=\mbox{\boldmath $0$}.
$$
So choosing
$$\displaystyle
q=2\mu\nabla\cdot\mbox{\boldmath $\phi$},
$$
we obtain the conclusion.

\noindent
$\Box$

See, when the case $\mbox{\boldmath $F$}=\mbox{\boldmath $0$}$, for example, \cite{FIIT}(1979), (8.7) and (8.8) on page 232 without mentioning any reference, and \cite{TCB} (1982) therein, some argument for the derivation of the solution together with a comparison of other arguments are given.  The solution in the case when $\mbox{\boldmath $F$}\equiv \mbox{\boldmath $0$}$ is a special case of the Papkovich-Neubar solution, see \cite{TCB}.
However, to my best knowledge, we cannot find the case when $\mbox{\boldmath $F$}\not\equiv \mbox{\boldmath $0$}$.  For example, in \cite{PRA} the authors seek a complete general solution of the Stokes system, that is,
in their words,  `other every solution can be obtained from it' under always the condition 
$\mbox{\boldmath $F$}\equiv \mbox{\boldmath $0$}$.  So maybe seeking the type of Proposition 2.1 is out of their scope.

However,  here we present how one can derive the Stokeslet by using Proposition 2.1.
\proclaim{\noindent Corollary 2.1.}
Let  $\mbox{\boldmath $\phi$}$ be an arbitrary  solution of
$$
\Delta\mbox{\boldmath $\phi$}-\mbox{\boldmath $a$}\,\delta(x)=\mbox{\boldmath $0$}.
\tag {2.3}
$$
If
$$\left\{
\begin{array}{l}
\displaystyle
\mbox{\boldmath $u$}=\frac{1}{\mu}
\left(\nabla(x\cdot\mbox{\boldmath $\phi$})-2\mbox{\boldmath $\phi$}\,\right),
\\
\\
\displaystyle
r=\nabla\cdot\mbox{\boldmath $\phi$},
\end{array}
\right.
\tag {2.4}
$$
then the pair $(\mbox{\boldmath $u$},r)$ satisfies (2.2).

\endproclaim
{\it\noindent Proof.}
In Proposition 2.1 choose $\mbox{\boldmath $F$}(x)=\mbox{\boldmath $a$}\,\delta(x)$.

\noindent
We have
$$\displaystyle
x\cdot\mbox{\boldmath $F$}
=(x\cdot\mbox{\boldmath $a$})\delta(x)=0
$$
and
$$\begin{array}{ll}
\displaystyle
(\nabla\mbox{\boldmath $F$})^Tx
&
\displaystyle
=(\mbox{\boldmath $a$}\otimes\nabla\delta(x))^Tx\\
\\
\displaystyle
&
\displaystyle
=(x\cdot\mbox{\boldmath $a$})\nabla\delta(x)\\
\\
\displaystyle
&
\displaystyle
=\nabla((x\cdot\mbox{\boldmath $a$})\delta(x))-\mbox{\boldmath $a$}\delta(x)\\
\\
\displaystyle
&
\displaystyle
=-\mbox{\boldmath $a$}\delta(x)\\
\\
\displaystyle
&
\displaystyle
=-\mbox{\boldmath $F$}.
\end{array}
$$
Thus the $(\mbox{\boldmath $v$},q)$ in Proposition 2.1 satisfies
$$\displaystyle
\mu\Delta\mbox{\boldmath $v$}-\nabla q+2\mu\mbox{\boldmath $a$}\delta(x)=\mbox{\boldmath $0$}.
$$
Divide the both sides by $2\mu$ and let
$$\left\{
\begin{array}{l}
\displaystyle
\mbox{\boldmath $u$}=\frac{\mbox{\boldmath $v$}}{2\mu},
\\
\\
\displaystyle
r=\frac{q}{2\mu}.
\end{array}
\right.
$$
Then the pair $(\mbox{\boldmath $u$},r)$ has the expression (2.4) and satisfies  (2.2).

\noindent
$\Box$

$\quad$

{\bf\noindent Remark 2.1.}
Needless to say, the author does not think Corollary 2.1 is not known at all.
The key of the proof above is the equations for the vector-valued  Dirac delta function 
$\mbox{\boldmath $F$}(x)=\mbox{\boldmath $a$}\,\delta(x)$:
$$\left\{\begin{array}{l}
\displaystyle
x\cdot\mbox{\boldmath $F$}=0,\\
\\
\displaystyle
(\nabla\mbox{\boldmath $F$})^Tx+\mbox{\boldmath $F$}=\mbox{\boldmath $0$}.
\end{array}
\right.
$$
As described in \cite{TCB} an origin of the solution type in Proposition 2.1 in the case when $\mbox{\boldmath $F$}\equiv\mbox{\boldmath $0$}$ goes back to the Papkovich-Neubar solution for the elasticity appeared in 1930s
(see the references therein).
In that time the Schwartz distribution theory yet appeared.  So the above equation did not have any 
exact meaning at that time.  Is this a reason
why one cannot find the statement like Corollary 2.1 together with 
Proposition 2.1? 
Or perhaps the appearance of the term $\mu\mbox{\boldmath $F$}-\mu(\nabla\mbox{\boldmath $F$})^Tx$
in the derived equation $(\mbox{\boldmath $v$},q)$ in Proposition 2.1 did not interest researchers in this field.

$\quad$

{\bf\noindent Example 1.}
Given  $\mbox{\boldmath $a$}$ a constant vector,
one can choose the vector field $\mbox{\boldmath $\phi$}$ satisfying (2.3) in Corollary 2.1 as
$$\displaystyle
\mbox{\boldmath $\phi$}(x)=-\frac{\mbox{\boldmath $a$}}{4\pi\vert x\vert}.
$$
Then the pair $(\mbox{\boldmath $u$},r)$ given by (2.4) takes the form
$$\begin{array}{ll}
\displaystyle
\mbox{\boldmath $u$}
&
\displaystyle
=-\frac{1}{8\pi\mu}\nabla\left(\frac{x\cdot\mbox{\boldmath $a$}}{\vert x\vert}\right)+\frac{2}{8\pi\mu}\frac{\mbox{\boldmath $a$}}{\vert x\vert}
\\
\\
\displaystyle
&
\displaystyle
=\frac{1}{8\pi\mu}\frac{\mbox{\boldmath $a$}}{\vert x\vert}
+\frac{1}{8\pi\mu}\frac{x\otimes x}{\vert x\vert}\mbox{\boldmath $a$}
\\
\\
\displaystyle
&
\displaystyle
=\mbox{\boldmath $J$}(x)\mbox{\boldmath $a$}
\end{array}
$$
and
$$\begin{array}{ll}
\displaystyle
r
&
\displaystyle
=\nabla\cdot\mbox{\boldmath $\phi$}
\\
\\
\displaystyle
&
\displaystyle
=
-\nabla\cdot\left(\frac{\mbox{\boldmath $a$}}{4\pi\vert x\vert}\right)
\\
\\
\displaystyle
&
\displaystyle
=\mbox{\boldmath $p$}(x)\cdot\mbox{\boldmath $a$}.
\end{array}
$$
So, by Corollary 2.1, we have reconfirmed that the pair $(\mbox{\boldmath $u$}, r)$ given by  (2.1) satisfies (2.2).

$\quad$

{\bf\noindent Example 2. Yarmukhamedov's Stokeslet.}
In 1972, Yarmukhamedov \cite{Y1} introduced a special fundamental solution of the Laplace equation
parametrized by a family of entire functions.
It takes the form
$$\displaystyle
-2\pi^2\,\Phi_K(x)=\int_0^{\infty}\,\text{Im}\,\left(\frac{K(w)}{w}\,\right)\frac{du}{\sqrt{\vert x'\vert^2+u^2}},
\tag {2.5}
$$
where $x=(x_1,x_2,x_3)$, $x'=(x_1,x_2)\not=(0,0)$, $w=x_3+i\sqrt{\vert x'\vert^2+u^2}$, $i=\sqrt{-1}$.

The function $K$ of the complex variable $w$ is an arbitrary entire function satisfying

$\bullet$  $\overline{K(w)}=K(\overline{w})$,

$\bullet$  $K(0)=1$,

$\bullet$  for each $R>0$ and $m=0,1,2$
$$\displaystyle
\sup_{\vert\text{Re}\,w\vert<R}\,\vert K^{(m)}(w)\vert<\infty.
$$

\noindent
Note that $K(w)\equiv 1$ satisfies all the conditions above and we have
$$\displaystyle
\Phi_K(x)=\frac{1}{4\pi\vert x\vert}.
$$

In an application to the Cauchy problem for the Laplace equation 
Yarmukhamedov in \cite{Y2} for a spacial $K$ and later, general $K$ in \cite{Y4}
had proven that there is a unique {\it entire harmonic function} $H_K(x)$ such that
$$\begin{array}{ll}
\displaystyle
\Phi_K(x)=\frac{1}{4\pi\vert x\vert}+H_K(x), 
& 
x'\not=(0,0).
\end{array}
\tag {2.6}
$$

$\quad$

{\bf\noindent Remark 2.2.}  This is a remark on a history.
Yarmukhamedov's original proof given in \cite{Y3} for general $K$ is incorrect.  This has been pointed out 
by the author in \cite{Islab} and its correct proof was given 
in Appendix A therein and in Appendix of \cite{IProbeCarleman} including also for the Helmholtz equation case.
In \cite{Islab} it is pointed out also that the $\Phi_K$ for specially chosen $K$ coincides with 
the Faddeev-Green function for the Laplace equation which played the central role in construction
of the complex geometrical optics solution \cite{SU2} in the study of  the Calder\'on problem \cite{Cal}.  In this sense, 
Yarmukhamedov's fundamental solution itself  is an interesting object.

$\quad$

The equation (2.6)  means that the `singularity' in the sense that integrand of $\Phi_K(x)$ in (2.6)
which is not absolutely integrable for $x'=(0,0)$, is removable
except for $x=(0,0,0)$
and the right hand side on (2.6) gives the unique extension of $\Phi_K$ to the whole space as the distribution.
Needless to say the origin $x=(0,0,0)$ is its unique singular point as a member of the singular support.
We denote the extension by the same symbol
$\Phi_K$.  Then $\Phi_K$ satisfies
$$\displaystyle
\Delta\Phi_K(x)+\delta(x)=0
$$
in the sense of distribution in the whole space.

Now let
$$\displaystyle
\mbox{\boldmath $\phi$}(x)=-\mbox{\boldmath $a$}\Phi_K(x).
$$
Then, the pair
$(\mbox{\boldmath $u$}, r)=(\mbox{\boldmath $u$}_K,
r_K)$ defined by (2.4) automatically satisfies (2.2).
Note that, by Example 1, we have the expression
$$\left\{
\begin{array}{l}
\displaystyle
\mbox{\boldmath $u$}_K(x)=\mbox{\boldmath $J$}(x)\mbox{\boldmath $a$}
-\frac{1}{\mu}
\left\{\nabla(x\cdot\mbox{\boldmath $a$}\,H_K(x))-2\mbox{\boldmath $a$}H_K(x)\,\right\},
\\
\\
\displaystyle
r_K(x)=\mbox{\boldmath $p$}(x)\cdot\mbox{\boldmath $a$}-
\nabla\cdot(\mbox{\boldmath $a$}\,H_K(x)).
\end{array}
\right.
\tag {2.7}
$$
It is clear that, by replacing $H_K$ on (2.7) with an arbitrary harmonic function, one gets another fundamental solution for the Stokes system.
However, the $H_K$ is important and interesting because it has a variety of applications, under various choice of $K$
together with a rotation of the original coordinates,
not only to the Cauchy problem for the Laplace equation \cite{Y4} and Helmholtz equation 
$\Delta u+k^2 u=0$ with $k>0$\footnote{Note that $\Delta=\frac{\partial^2}{\partial x_1^2}
+\frac{\partial^2}{\partial x_2^2}+\frac{\partial^2}{\partial x_3^2}$ and thus
this is not the modified Helmholtz equation
$\Delta u-k^2u=0$.}
(see Remark 5.4 in \cite{IProbeCarleman}), but also the enclosure method for
inverse conductivity problem in the infinte slabe \cite{Islab} (see also Subsection 3.1.4 of \cite{IReview}),
and the probe method itself \cite{IProbeCarleman} described in Subsection 2.4.

It should be noted that, as a direct and simple application of (2.7) with a specially chosen $K$
(see (2.11) for a choice), one can obtain an extension of a formula of the Carleman-type for the Laplace equation \cite{Y4} to the Stokes system.
Since we are concentrated on the probe and singular sources methods themselves, we do not mention about it any more.
Instead, for a study in two dimensions with a different approach, see \cite{Ar}.

\subsection{Needle sequence}

Here we introduce the key concept of  the probe method reformulated in \cite{INew}.

$\quad$

{\bf\noindent Definition 2.1.}  Fix $\mbox{\boldmath $a$}\not=\mbox{\boldmath $0$}$.
Let $x$ and $\sigma\in N_x$.  We call the sequence $\{(\mbox{\boldmath $v$}_n, q_n)\}$ of  $H^2(\Omega)\times H^1(\Omega)$ solutions
of the Stokes system
a needle sequence for $(x,\sigma)$  with respect to the singular solution given by $(\mbox{\boldmath $J$}(\,\cdot\,-x)\mbox{\boldmath $a$}, \mbox{\boldmath $p$}(\,\cdot\,-x)\cdot\mbox{\boldmath $a$}))$,
 if it satisfies 
for each nonempty open set $U$ of  $\Bbb R^3$ with $\overline{U}\subset\Omega\setminus\sigma$
$$\displaystyle
\lim_{n\rightarrow\infty}\Vert\mbox{\boldmath $v$}_n-\mbox{\boldmath $J$}(\,\cdot\,-x)\mbox{\boldmath $a$}\Vert_{H^2(U)}
+\Vert q_n-\mbox{\boldmath $p$}(\,\cdot\,-x)\cdot\mbox{\boldmath $a$})\Vert_{H^1(U)}=0.
\tag {2.8}
$$
This means that the sequence $\{(\mbox{\boldmath $v$}_n, q_n)\}$ converges to $(\mbox{\boldmath $J$}(\,\cdot\,-x)\mbox{\boldmath $a$}, \mbox{\boldmath $p$}(\,\cdot\,-x)\cdot\mbox{\boldmath $a$})$
in $H^2_{\text{loc}}(\Omega\setminus\sigma)\times H^1_{\text{loc}}(\Omega\setminus\sigma)$.
We denote by ${\cal N}(x,\sigma)={\cal N}(x,\sigma,\mbox{\boldmath $a$})$ the set of all needle sequences for $(x,\sigma)$.

$\quad$

In this paper, for a fixed $\mbox{\boldmath $a$}$ the $\mbox{\boldmath $J$}(\,\cdot\,-x)\mbox{\boldmath $a$}$ plays the central role.
So we start with describing the following simple, however, important fact.

\proclaim{\noindent Proposition 2.2.}  Fix $\mbox{\boldmath $a$}\not=\mbox{\boldmath $0$}$.
We have, for any finite cone $V$ with vertex at  $x$ 
$$\displaystyle
\int_{V}\vert\text{Sym}\,\nabla(\mbox{\boldmath $J$}(y-x)\mbox{\boldmath $a$})\vert^2\,dy=\infty.
$$

\endproclaim

{\it\noindent Proof.}
Without loosing a generality, one may consider only $x=(0,0,0)$.
A direct computation yields
$$\displaystyle
\nabla\left(\frac{\mbox{\boldmath $a$}}{\vert y\vert}\,\right)=-\frac{\mbox{\boldmath $a$}\otimes y}{\vert y\vert^3}
$$
and
$$
\displaystyle
\nabla\left(\frac{y\otimes y}{\vert y\vert^3}\,\mbox{\boldmath $a$}\,\right)
=\frac{y\cdot\mbox{\boldmath $a$}}{\vert y\vert^3}\,I_3+\frac{y\otimes\mbox{\boldmath $a$}}{\vert y\vert^3}
-3(y\cdot\mbox{\boldmath $a$})\frac{y\otimes y}{\vert y\vert^5}.
$$
Thus one gets
$$\displaystyle
8\pi\mu\,\nabla(\mbox{\boldmath $J$}(y)\mbox{\boldmath $a$})
=-\frac{\mbox{\boldmath $a$}\otimes y}{\vert y\vert^3}+\frac{y\cdot\mbox{\boldmath $a$}}{\vert y\vert^3}\,I_3+\frac{y\otimes\mbox{\boldmath $a$}}{\vert y\vert^3}
-3(y\cdot\mbox{\boldmath $a$})\frac{y\otimes y}{\vert y\vert^5}
$$
and 
$$\displaystyle
8\pi\mu\,\text{Sym}\,(\nabla\mbox{\boldmath $J$}(y)\mbox{\boldmath $a$})
=\frac{y\cdot\mbox{\boldmath $a$}}{\vert y\vert^3}\,I_3-3(y\cdot\mbox{\boldmath $a$})\frac{y\otimes y}{\vert y\vert^5}.
$$
Therefore we have the expression
$$\displaystyle
(8\pi\mu)^2\vert\,\text{Sym}\,(\nabla\mbox{\boldmath $J$}(y)\mbox{\boldmath $a$})\vert^2
=\frac{4(y\cdot\mbox{\boldmath $a$})^2}{\vert y\vert^6}.
$$
Here, the set $V\cap S^2$ is open and nonempty. Then, for $\frac{\mbox{\boldmath $a$}}{\vert\mbox{\boldmath $a$}\vert}\in S^2$ one can find
a nonempty open subset $U$ of  $V\cap S^2$ such that $\inf_{\omega\in U}\vert\mbox{\boldmath $a$}\cdot\omega\vert>0$.
This yields
$$\displaystyle
0<\int_{V\cap S^2}(\mbox{\boldmath $a$}\cdot\mbox{\boldmath $\omega$})^2 d\mbox{\boldmath $\omega$}.
$$
Thus one gets
$$\displaystyle
\int_{V}\frac{(y\cdot\mbox{\boldmath $a$})^2}{\vert y\vert^6}\,dy
=\int_{0}^{R}\frac{dr}{r^2}\int_{V\cap S^2}(\mbox{\boldmath $a$}\cdot\mbox{\boldmath $\omega$})^2 d\mbox{\boldmath $\omega$}
=\infty.
$$

\noindent
$\Box$

As well as Lemmas 2.1-2.2 of \cite{INew} which are the Helmholtz equation case, 
it is important to check the two types of blowing up property described below.

\proclaim{\noindent Theorem 2.1.}  Fix $\mbox{\boldmath $a$}\not=\mbox{\boldmath $0$}$.
Given an arbitrary point $x\in\Omega$ and needle $\sigma\in N_x$ let $\{(\mbox{\boldmath $v$}_n, q_n)\}\in {\cal N}(x,\sigma)$.

\noindent
(a)  Let $V$ be an arbitrary finite cone with vertex at $x$.  Then, we have
$$\displaystyle
\lim_{n\rightarrow\infty}\Vert \text{Sym}\,\nabla\mbox{\boldmath $v$}_n\Vert_{L^2(V\cap\Omega)}=\infty.
$$

\noindent
(b)  Let $z\in\Omega$ be an arbitrary point on $\sigma\setminus\{x\}$ and $B$ an open ball centered at $z$.
Then, we have
$$\displaystyle
\lim_{n\rightarrow\infty}\Vert \text{Sym}\,\nabla\mbox{\boldmath $v$}_n\Vert_{L^2(B\cap\Omega)}=\infty.
$$

\endproclaim

The property (a) compared with (b) tells us a more precise blowing up of the needle sequence at the tip of the needle.
This needs to treat the caee when only the tip of the needle is located on the boundary of an obstacle.

Except for the proof of  (a) which is a corollary of Proposition 2.2 with the help of Fatou's lemma,
unlike the Helmholtz equation case,  the proof of (b) needs a special care 
caused by the system, that is, the appearance of the pressure term $q_n$ and as a result we take a different way from the original proof of Lemma 2.2 in \cite{INew}.
The proof is given in Subsection 4.1.
Note that the behaviour of $q_n$ on $\sigma$ is not important, rather, we never use it.

As a corollary of Theorem 2.1 and Definition 2.1 
we have the characterization of  needle $\sigma$ by using an
arbitrary fixed needle sequence $\{(\mbox{\boldmath $v$}_n,q_n)\}$ for $(x,\sigma)$:
$$\displaystyle
\sigma(]0,\,1])=
\left\{\,z\in\Omega\,\vert\,\forall B\in\beta_z
\,
\lim_{n\rightarrow\infty}\Vert \text{Sym}\,\nabla\mbox{\boldmath $v$}_n\Vert_{L^2(B\cap\Omega)}=\infty\,
\right\},
\tag {2.9}
$$
where $\beta_z$ denotes the set of all open balls centered at $z$.
This means that the needle sequence does not loose the information of the geometry of the original needle.
The needle sequence $\{(\mbox{\boldmath $v$}_n,q_n)\}$ is a kind of the (real-valued)
solutions of the Stokes system 
whose energies are concentrated on $\sigma$.  The author thinks such solution itself  has independent interest.

\subsection{Existence of the needle sequence}

The needle sequence {\it exists}.  It is a corollary of Proposition 2.1 which makes us possible to use the needle sequence for the Laplace equation.  
Let us explain it here.

Given $x\in\Omega$ and $\sigma\in N_x$
let $\{v_n\}$ be a needle sequence for $(x,\sigma)$ for the Laplace equation \cite{INew}, that satisfies,
for each nonempty open set $U$ of  $\Bbb R^3$ with $\overline{U}\subset\Omega\setminus\sigma$
$$\displaystyle
\lim_{n\rightarrow\infty}\left\Vert v_n-\frac{1}{4\pi\vert \,\cdot\,-x\vert}\right\Vert_{H^1(U)}=0.
$$
The existence of this sequence is a consequence of the Runge approximation property of the Laplace equation.
Note that by this convergence, one can replace the space $H^1(U)$ with $H^m(U)$ for all $m\ge 1$ and, in particular,
$m=3$.  

Substituting $-\mbox{\boldmath $a$}\,v_n$ into $\mbox{\boldmath $\phi$}$ in (2.4) of Corollary 2.1,
we see that the sequence $\{\mbox{\boldmath $v$}_n,q_n)\}$
given by
$$\left\{\begin{array}{l}
\displaystyle
\mbox{\boldmath $v$}_n(y)
=-\frac{1}{\mu}\{\nabla(y\cdot\mbox{\boldmath $a$}\,v_n(y))-2\mbox{\boldmath $a$}v_n(y)\},
\\
\\
\displaystyle
q_n(y)=-\nabla\cdot(\mbox{\boldmath $a$}\,v_n(y)),
\end{array}
\right.
\tag {2.10}
$$
satisfies the convergence property (2.8) for each nonempty open set $U$
of $\Bbb R^3$ with $\overline{U}\subset\Omega\setminus\sigma$.
Thus $\{(\mbox{\boldmath $v$}_n, q_n)\}\in{\cal N}(x,\sigma)$.

\subsection{Concrete example of the needle sequence}

The concrete example of the needle sequence for the Laplace equation,
for a simple, however, important needle has been given.
For simplicity of description and without loosing generality, one may assume $x=(0,0,0)\in\Omega$ and
 $\sigma\in N_x$ having the form
 $$\displaystyle
 \sigma=l\cap\Omega
 $$
 where $l=\{y\in\Bbb R^3\,\vert\,y=s\mbox{\boldmath $e$}_3,\, s\ge 0\}$.
 We assume that $l\cap\Omega$ is connected.  This type of needle with tip at $x=(0,0,0)$ is called the straight needle
 directed $\mbox{\boldmath $e$}_3$.
 
 Let $E_{\alpha}(w)$ denote the Mittag-Leffler function with parameter $\alpha$ restricted $0<\alpha<1$:
 $$\displaystyle
 E_{\alpha}(w)=\sum_{m=0}^{\infty}\frac{w^m}{\Gamma(\alpha+1)}.
 $$
 Define
 $$\displaystyle
 K(w)=E_{\alpha}(\tau w),
 \tag {2.11}
 $$
 where $\tau$ is a positive number.  
 The $K$ defined by (2.11) is an entire function and satisfies the three conditions listed in Example 2.
 Thus one can construct $\Phi_K$ by substituting this $K$ into (2.5) and we have
 the harmonic function $H_K$ in the whole space.  
 
 In \cite{IProbeCarleman} as a corollary
 of Theorem 4.2 therein
 the author probed that
 by choosing a suitable sequence $\alpha=\alpha_n\rightarrow 0$ and $\tau=\tau_n\rightarrow\infty$, $n=1,2,\cdots$
 the sequence
 $$\begin{array}{ll}
 \displaystyle
 v_n(y)=-H_{K_n}(y), & 
 \displaystyle
 n=1,2,\cdots,
 \end{array}
 \tag {2.12}
 $$
 yields a needle sequence for $(x,\sigma)$ for the Laplace equation.  Here $K_n(w)=E_{\alpha_n}(\tau_nw)$.
 Therefore, by substituting this $v_n$ into the formula (2.10), one gets a concrete sequence 
 $\{(\mbox{\boldmath $v$}_n,q_n)\}\in{\cal N}(x,\sigma)$.

Finally we note that the $v_n$ given by (2.12) has the explicit integral representation \cite{IProbeCarleman}, (4.11) on page 1881:
$$\displaystyle
v_n(y)
=\frac{1}{2\pi^2}\,\int_0^{\infty}
\text{Im}\,\left(\,\frac{E_{\alpha_n}(\tau_n w)-1}{w}\,\right)\,\frac{du}{\sqrt{\vert y'\vert^2+u^2}},
$$
where $w=y_3+i\sqrt{\vert y'\vert^2+u^2}$ and thus by (2.10), we can compute 
each $(\mbox{\boldmath $v$}_n, q_n)$ constructed above  in principle.  
Besides, making a rotation and moving in parallel,
one gets an explicit needle sequence for the straight needle directed to an arbitrary direction and with
an arbitrary point as the tip.

\section{Integrated theory}

First, we introduce a vector-valued function that plays a central role in both methods.
For the purpose we note that 
 $$\begin{array}{ll}
 \displaystyle
 \int_{\partial D}\mbox{\boldmath $J$}(z-x)\mbox{\boldmath $a$}\cdot\mbox{\boldmath $\nu$}\,dS=0,
 & x\in\Omega\setminus\overline{D}.
 \end{array}
 $$
 since we have $\nabla\cdot \mbox{\boldmath $J$}(\,\cdot\,-x)\mbox{\boldmath $a$}=0$ in $D$.

$\quad$

{\bf\noindent Definition 3.1.}
Let $\mbox{\boldmath $a$}$ be an arbitrary constant vector.
Given $x\in\Omega\setminus\overline{D}$ let the pair $(\mbox{\boldmath $w$},q)\in H^2(\Omega\setminus\overline{D})\times H^1(\Omega\setminus\overline{D})$ be
the solution of
$$\left\{
 \begin{array}{ll}
 \mu\Delta\mbox{\boldmath $w$}-\nabla q=\mbox{\boldmath $0$}, & z\in\Omega\setminus\overline{D},
 \\
 \\
 \displaystyle
 \nabla\cdot\mbox{\boldmath $w$}=0, &  z\in\Omega\setminus\overline{D},\\
 \\
 \displaystyle
 \mbox{\boldmath $w$}=-\mbox{\boldmath $J$}(z-x)\mbox{\boldmath $a$}, & z\in\partial D,
 \\
 \\
 \displaystyle
 \mbox{\boldmath $w$}=\mbox{\boldmath $0$}, & z\in\partial\Omega.
 \end{array}
 \right.
 \tag {3.1}
 $$
The $\mbox{\boldmath $w$}$ is unique, depends on $x$ and linear with respect to $\mbox{\boldmath $a$}$.
Thus, one has the expression
 $$\displaystyle
 \mbox{\boldmath $w$}=\mbox{\boldmath $w$}_x\mbox{\boldmath $a$},
 $$
 where $\mbox{\boldmath $w$}_x$ is a $M_3(\Bbb R)$-valued function.  We call the vector-valued function
 $\mbox{\boldmath $w$}_x\mbox{\boldmath $a$}$ the reflected solution by $D$ caused by the incident incident filed\footnote{
 Later we consider also this incident field replaced with another one.} $\mbox{\boldmath $J$}(\,\cdot\,-x)\mbox{\boldmath $a$}$.
 The $q$ is unique modulo a constant function.  We denote it by the symbol $q_x=q_x(\,\cdot\,;\mbox{\boldmath $a$})$.

\subsection{The probe method}

In this section we briefly describe the probe method for the Stokes system.

\subsubsection{Side A}

The most fundamental concept in the probe method is defined as follows.

$\quad$

{\bf\noindent Definition 3.1.}
Fix  $\mbox{\boldmath $a$}\not=\mbox{\boldmath $0$}$.
Given $x\in\Omega$, $\sigma\in N_x$ and
$\mbox{\boldmath $\xi$}=\{(\mbox{\boldmath $v$}_n, q_n)\}
\in {\cal N}(x,\sigma,\mbox{\boldmath $a$})$,
define the sequence given by
$$\displaystyle
I(x,\sigma, \mbox{\boldmath $a$},\mbox{\boldmath $\xi$})_n=
<(\Lambda_0-\Lambda_D)
\mbox{\boldmath $v$}_n\vert_{\partial\Omega},
\mbox{\boldmath $v$}_n\vert_{\partial\Omega}>,
$$
where this right-hand side has the meaning given by (1.12).
We call this the {\it indicator sequence} for the probe method.

$\quad$

From (1.10) we have the expression
$$\displaystyle
I(x,\sigma, \mbox{\boldmath $a$},\mbox{\boldmath $\xi$})_n
=-\Vert 2\mu\,\text{Sym}\nabla \mbox{\boldmath $w$}_n\Vert_{L^2(\Omega\setminus\overline D)}
-\Vert 2\mu\,\text{Sym}\,\nabla\mbox{\boldmath $v$}_n(\,\cdot\,;\sigma,\mbox{\boldmath $a$}\Vert_{L^2(D)}^2,
\tag {3.2}
$$
where $\mbox{\boldmath $w$}_n=\mbox{\boldmath $u$}_n-\mbox{\boldmath $v$}_n$,
$\mbox{\boldmath $v$}_n=\mbox{\boldmath $v$}_n(\,\cdot\,;\sigma,\mbox{\boldmath $a$})$ and $\mbox{\boldmath $u$}_n$
and a $p=p_n$ solve (1.6) with $\mbox{\boldmath $f$}=\mbox{\boldmath $v$}_n$ on $\partial\Omega$.

The Side A of the probe method is concerned with the convergence property
of the indicator sequence together with the behaviour of  a function appeared as the limit.
It is based on the following fact.

\proclaim{\noindent Proposition 3.1.}
Given $x\in\Omega\setminus\overline{D}$ and $\sigma\in N_x$, let $\{(\mbox{\boldmath $v$}_n,q_n)\}\in {\cal N}(x,\sigma)$.
Let $(\mbox{\boldmath $u$},p)=(\mbox{\boldmath $u$}_n,p_n)$ solve (1.1) with $\mbox{\boldmath $f$}=\mbox{\boldmath $v$}_n$
on $\partial\Omega$.
If $\sigma\cap\overline{D}=\emptyset$, then we have, in the $H^2(\Omega\setminus\overline{D})$-topology
$$\displaystyle
\lim_{n\rightarrow\infty}(\mbox{\boldmath $u$}_n-\mbox{\boldmath $v$}_n)=\mbox{\boldmath $w$}_x\mbox{\boldmath $a$}.
$$
Besides, for all $\mbox{\boldmath $g$}\in H^{\frac{1}{2}}(\partial\Omega)$ with $\int_{\partial\Omega}\mbox{\boldmath $g$}\cdot\mbox{\boldmath $\nu$}\,dS=0$
we have
$$\displaystyle
\lim_{n\rightarrow\infty}<(\Lambda_0-\Lambda_D)\mbox{\boldmath $v$}_n\vert_{\partial\Omega},\mbox{\boldmath $g$}>
=-\int_{\partial\Omega}\sigma(\mbox{\boldmath $w$}_x(z)\mbox{\boldmath $a$}, q_x(z;\mbox{\boldmath $a$}))\mbox{\boldmath $\nu$}\cdot\mbox{\boldmath $g$}\,dS(z).
\tag {3.3}
$$

\endproclaim

{\it\noindent Proof.}
The function $\mbox{\boldmath $w$}_n=\mbox{\boldmath $u$}_n-\mbox{\boldmath $v$}_n$ satisfies
$$\left\{
\begin{array}{ll}
\displaystyle
\mu\Delta\mbox{\boldmath $w$}_n-\nabla (p_n-q_n)=\mbox{\boldmath $0$}, & y\in\Omega\setminus\overline{D},
\\
\\
\displaystyle
\nabla\cdot\mbox{\boldmath $w$}_n=0, & y\in\Omega\setminus\overline{D},\\
\\
\displaystyle
\mbox{\boldmath $w$}_n=-\mbox{\boldmath $v$}_n, & y\in\partial D,
\\
\\
\displaystyle
\mbox{\boldmath $w$}_n=\mbox{\boldmath $0$}, & y\in\partial\Omega.
\end{array}
\right.
$$
Since $\mbox{\boldmath $v$}_n$ converges to $\mbox{\boldmath $J$}(\,\cdot-x)\mbox{\boldmath $a$}$ in $H^2(D)$
as $n\rightarrow\infty$,
we see that $\mbox{\boldmath $w$}_n$ and $p_n-q_n-(p_n-q_n)_{E}$ with $E=\Omega\setminus\overline{D}$ are
Cauchy sequences in $H^2(\Omega\setminus\overline{D})$ and $H^1(\Omega\setminus\overline{D})$, respectively.
Thus $\mbox{\boldmath $w$}_n$ converges to a $\mbox{\boldmath$w$}_{\infty}\in H^2(\Omega\setminus\overline{D})$ and $p_n-q_n-(p_n-q_n)_{E}\rightarrow r\in H^1(\Omega\setminus\overline{D})$.  Therefore one gets
$$\left\{
\begin{array}{ll}
\displaystyle
\Delta\mbox{\boldmath $w$}_{\infty}-\nabla r=\mbox{\boldmath $0$}, & y\in\Omega\setminus\overline{D},
\\
\\
\displaystyle
\nabla\cdot\mbox{\boldmath $w$}_{\infty}=\mbox{\boldmath $0$}, & y\in\Omega\setminus\overline{D},\\
\\
\displaystyle
\mbox{\boldmath $w$}_{\infty}=-\mbox{\boldmath $J$}(y-x)\mbox{\boldmath $a$} & y\in\partial D,
\\
\\
\displaystyle
\mbox{\boldmath $w$}_{\infty}=\mbox{\boldmath $0$}, & y\in\partial\Omega.
\end{array}
\right.
$$
By the uniqueness of the velocity field we obtain $\mbox{\boldmath $w$}_{\infty}=\mbox{\boldmath $w$}_x\mbox{\boldmath $a$}$.
Thus the first assertion is valid.
And also $r=q_x+C$ with a constant $C$.
Then we have
$$\begin{array}{l}
\displaystyle
\,\,\,\,\,\,
\lim_{n\rightarrow\infty}<(\Lambda_0-\Lambda_D)\mbox{\boldmath $v$}_n\vert_{\partial\Omega},
\mbox{\boldmath $g$}>
\\
\\
\displaystyle
=
-\lim_{n\rightarrow\infty}\int_{\partial\Omega}\sigma(\mbox{\boldmath $u$}_n-\mbox{\boldmath $v$}_n, p_n-q_n)\mbox{\boldmath $\nu$}\cdot\mbox{\boldmath $g$}\,dS(z)
\\
\\
\displaystyle
=-\int_{\partial\Omega}\sigma(\mbox{\boldmath $w$}_x(z)\mbox{\boldmath $a$}, q_x(z)+C)
\mbox{\boldmath $\nu$}\cdot\mbox{\boldmath $g$}\,dS(z).
\end{array}
$$
This together with the condition $\int_{\partial\Omega}\mbox{\boldmath $g$}\cdot\mbox{\boldmath $\nu$}\,dS=0$ yields the validity of the second assertion.

\noindent
$\Box$

The first fact in Proposition 3.1 together with identity (1.10) {\it suggests} us to introduce the following function.

$\quad$

{\bf\noindent Definition 3.2.}  Let $\mbox{\boldmath $a$}$ and $\mbox{\boldmath $b$}$ be arbitrary constant vectors.
Given $(x,y)\in(\Omega\setminus\overline{D})^2$ define
$$\begin{array}{ll}
\displaystyle
I(x,\mbox{\boldmath $a$};y,\mbox{\boldmath $b$})
&
\displaystyle
=-\int_{\Omega\setminus\overline{D}}
2\mu\,\text{Sym}\,\nabla\mbox{\boldmath $w$}_x\mbox{\boldmath $a$}
\cdot
\text{Sym}\,\nabla\mbox{\boldmath $w$}_y\mbox{\boldmath $b$}\,dz
\\
\\
\displaystyle
&
\displaystyle
\,\,\,
-\int_D\,2\mu\,\text{Sym}\nabla\mbox{\boldmath $J$}(z-x)\mbox{\boldmath $a$}\cdot\text{Sym}\,\nabla\mbox{\boldmath $J$}(z-y)\mbox{\boldmath $b$}\,dz
\end{array}
$$
and
$$\begin{array}{ll}
\displaystyle
I(x,\mbox{\boldmath $a$})
&
\displaystyle
\equiv I(x,\mbox{\boldmath $a$};x,\mbox{\boldmath $a$})
\\
\\
\displaystyle
&
\displaystyle
=-\Vert 2\mu\,\text{Sym}\,\nabla\mbox{\boldmath $w$}_x\mbox{\boldmath $a$}\Vert_{L^2(\Omega\setminus\overline{D})}^2
-\Vert 2\mu\,\text{Sym}\nabla\mbox{\boldmath $J$}(\,\cdot\,-x)\mbox{\boldmath $a$}\Vert_{L^2(D)}^2.
\end{array}
\tag {3.4}
$$
We call this the indicator function for the probe method and the function $I(x,\mbox{\boldmath $a$};y,\mbox{\boldmath $b$})$ the {\it lifting} of the indicator function.

$\quad$

Note that we have the trivial symmetry relation
$$\displaystyle
I(x,\mbox{\boldmath $a$};y,\mbox{\boldmath $b$})
=I(y,\mbox{\boldmath $b$};x,\mbox{\boldmath $a$}).
$$

The validity of  the following fact is now clear.

\proclaim{\noindent Proposition 3.2.}
Let $(x,y)\in(\Omega\setminus\overline{D})^2$ and $(\sigma,\sigma')\in N_x\times N_y$.
Let $\{(\mbox{\boldmath $v$}_n, q_n)\}\in {\cal N}(x,\sigma,\mbox{\boldmath $a$})$
and $\{(\mbox{\boldmath $v$}_n', q_n')\}\in{\cal N}(y,\sigma',\mbox{\boldmath $b$})$.
If $\sigma\cap\overline{D}=\emptyset$ and $\sigma'\cap\overline{D}=\emptyset$, then we have
$$\begin{array}{ll}
\displaystyle
I(x,\mbox{\boldmath $a$};y,\mbox{\boldmath $b$})
&
=
\lim_{n\rightarrow\infty}<(\Lambda_0-\Lambda_D)\mbox{\boldmath $v$}_n\vert_{\partial\Omega},\mbox{\boldmath $v$}_n'\vert_{\partial\Omega}>.
\end{array}
$$

\endproclaim

\noindent
This proposition means that the lifting of the indicator 
function $I(x,\mbox{\boldmath $a$};y,\mbox{\boldmath $b$})$ 
of $(x,y)\in (\Omega\setminus\overline{D})^2$ can be computed from
the data $\Lambda_0-\Lambda_D$, in principle, by using the needle sequences in
${\cal N}(x,\sigma,\mbox{\boldmath $a$})$ and $ {\cal N}(y,\sigma',\mbox{\boldmath $b$})$,
provided $\sigma\cap\overline{D}=\sigma'\cap\overline{D}=\emptyset$.

The theoretical core of the Side A of the probe method is summarized as follows.

\proclaim{\noindent Theorem 3.1.}  Fix $\mbox{\boldmath $a$}\not=\mbox{\boldmath $0$}$.

\noindent
(i)  Let $x\in\Omega\setminus\overline{D}$ and $\sigma\in N_x$ satisfy $\sigma\cap\overline{D}=\emptyset$.
Then, for all $\mbox{\boldmath $\xi$}\in {\cal N}(x,\sigma,\mbox{\boldmath $a$})$ we have
$$\displaystyle
\lim_{n\rightarrow\infty}\,I(x,\sigma, \mbox{\boldmath $a$},\mbox{\boldmath $\xi$})_n=
I(x,\mbox{\boldmath $a$}).
\tag {3.5}
$$

\noindent
(ii) We have
$$\displaystyle
\lim_{x\rightarrow a\in\partial D}I(x,\mbox{\boldmath $a$})=-\infty.
\tag {3.6}
$$

\noindent
(iii)  For each positive number $\epsilon$ it holds that
$$\displaystyle
-\infty<\inf_{x\in\Omega\setminus\overline{D},\,\text{dist}\,(x,\partial D)>\epsilon}\,I(x,\mbox{\boldmath $a$})<0.
\tag {3.7}
$$

\endproclaim

{\it\noindent Proof.}  (i)  is the particular case of Proposition 3.2.
For (ii), from (3.4) we have
$$\displaystyle
I(x,\mbox{\boldmath $a$})\le
-\Vert 2\mu\,\text{Sym}\,\nabla(\mbox{\boldmath $J$}(z-x)\mbox{\boldmath $a$})\Vert_{L^2(D)}^2.
\tag {3.8}
$$
A combination of (3.8) and Proposition 2.2 with an appropriate choice of a finite cone $V$ yields (ii).
Finally, as a consequence of the well-posedness of the Stokes system, we have
the boundedness of the $H^1(\Omega\setminus\overline{D})$-norm
of  reflected solution $\mbox{\boldmath $w$}_x\mbox{\boldmath $a$}$ with respect to $x\in\Omega\setminus\overline{D}$ with $\text{dist}(x,\partial D)>\epsilon$. 
Thus the expression (3.4) yields (iii).

\noindent
$\Box$

\noindent
Theorem 3.1 means that:
if $\sigma\cap\overline{D}=\emptyset$,
as $n\rightarrow\infty$
the indicator sequence $\{I(x,\sigma,\mbox{\boldmath $a$},\mbox{\boldmath $\xi$})_n\}$ is convergent
and its limit gives us the value of the indicator function, which blows up when the tip of $\sigma$, that is,
$x$, approaches a point on $\partial D$.

\noindent
This is the Side A of the probe method.

\subsubsection{Side B}

The Side B asks us to describe the asymptotic behaviour of the indicator sequence as $n\rightarrow$
without the condition $\sigma\cap\overline{D}=\emptyset$.
This consideration was not present in previous singular sources method, and in this respect the probe method has an advantage.

\proclaim{\noindent Theorem 3.2.}  Fix $\mbox{\boldmath $a$}\not=\mbox{\boldmath $0$}$.
Let $x\in\Omega$ and $\sigma\in N_x$.
Assume that one of two cases (a) and (b) listed below are satified:

\noindent
(a)  $x\in\overline{D}$;

 \noindent
 (b)  $x\in\Omega\setminus\overline{D}$ and $\sigma$ satisfies $\sigma\cap D\not=\emptyset$.
 
\noindent
Then, for any needle sequence $\mbox{\boldmath $\xi$}\in{\cal N}(x,\sigma,\mbox{\boldmath $a$})$ we have
$$\displaystyle
\lim_{n\rightarrow\infty}\,I(x,\sigma,\mbox{\boldmath $a$},\mbox{\boldmath $\xi$})_n=-\infty.
\tag {3.9}
$$
\endproclaim

{\it\noindent Proof.}
From (3.2) one gets
$$\displaystyle
I(x,\sigma,\mbox{\boldmath $a$},\mbox{\boldmath $\xi$})_n
\le
-\Vert 2\mu\,\text{Sym}\,\nabla\mbox{\boldmath $v$}_n(\,\cdot\,;\sigma,\mbox{\boldmath $a$})\Vert_{L^2(D)}^2.
$$
This together with the  blowing up property of the needle sequence in Theorem 2.1 yields (3.9).
Note that in particular, to treat the case $x\in\partial D$ we need (a) of  Theorem 2.1.

\noindent
$\Box$

\noindent
This is the Side B of the probe method. 
\footnote{It should be noted that until now, the {\it grazing case}, that is the case when all of the conditions
$x\in\Omega\setminus\overline{D}$, $\sigma\cap D=\emptyset$ and $\sigma\cap\partial D\not=\emptyset$ 
are satisfied, we do not know the behaviour of  the indicator sequence.}
 
As a corollary with the help pf the {\it existence} of the needle sequence,
using the same argument as \cite{IPS},
we have the characterization of $\overline{D}$ {\it in terms of} the blowing up property of the indicator sequence.
That is:

\proclaim{\noindent Corollary 3.1.}  Fix $\mbox{\boldmath $a$}\not=\mbox{\boldmath $0$}$.
A point $x\in\Omega$ belongs to $\overline{D}$
if and only if for any needle $\sigma\in N_x$ and needle 
sequence $\mbox{\boldmath $\xi$}\in{\cal N}(x,\sigma,\mbox{\boldmath $a$})$ the
indicator sequence $\{I(x,\sigma,\mbox{\boldmath $a$},\mbox{\boldmath $\xi$})_n\}$
blows up as $n\rightarrow\infty$.

\endproclaim

{\bf\noindent Remark 3.1.}
Note that the statement $x\in\overline{D}$ is independent of $\mbox{\boldmath $a$}$.
Thus Corollary 3.1 yields:
if, for a nonzero vector $\mbox{\boldmath $a$}_0$ we have for any needle $\sigma\in N_x$ and needle 
sequence $\mbox{\boldmath $\xi$}\in{\cal N}(x,\sigma,\mbox{\boldmath $a$}_0)$
indicator sequence $\{I(x,\sigma,\mbox{\boldmath $a$}_0,\mbox{\boldmath $\xi$})_n\}$
blows up as $n\rightarrow\infty$,
then for any nonzero vector $\mbox{\boldmath $a$}$, any needle $\sigma\in N_x$ and needle 
sequence $\mbox{\boldmath $\xi$}\in{\cal N}(x,\sigma,\mbox{\boldmath $a$})$
indicator sequence $\{I(x,\sigma,\mbox{\boldmath $a$},\mbox{\boldmath $\xi$})_n\}$
blows up as $n\rightarrow\infty$.
Therefore, the choice of $\mbox{\boldmath $a$}$ is no important.

$\quad$

Let $\mbox{\boldmath $a$}(\not=\mbox{\boldmath $0$})$ be an arbitrary fixed constant vector.
Corollary 3.1 yields the following  reconstruction formulae\footnote{Needless to say, this automatically yields the uniqueness of recovering of $D$ from $\Lambda_0-\Lambda_D$
since we have $D=\{x\in\overline{D}\,\vert\, \exists B\in\beta_x\,\, B\subset\overline{D}\}$.}
of $\overline{D}$:
$$\begin{array}{ll}
\displaystyle
\overline{D}
&
\displaystyle
=\{x\in\Omega\,\vert\,\forall\sigma\in N_x\,\forall\mbox{\boldmath $\xi$}\in{\cal N}(x,\sigma,\mbox{\boldmath $a$})\,
\lim_{n\rightarrow\infty}\,I(x,\sigma,\mbox{\boldmath $a$},\mbox{\boldmath $\xi$})_n=-\infty\,\}
\\
\\
\displaystyle
&
\displaystyle
=\{x\in\Omega\,\vert\,\forall\sigma\in N_x\,\forall\mbox{\boldmath $\xi$}\in{\cal N}(x,\sigma,\mbox{\boldmath $a$})\,
\lim_{n\rightarrow\infty}\,\vert I(x,\sigma,\mbox{\boldmath $a$},\mbox{\boldmath $\xi$})_n\vert=\infty\,\}.
\end{array}
$$

So the set $\overline{D}$ is transformed into the following set of numerical values without loosing the original geometric information:
$$\displaystyle
\left\{
\{I(x,\sigma,\mbox{\boldmath $a$},\mbox{\boldmath $\xi$})_n\}
\,\vert\,x\in\Omega,\sigma\in N_x, \mbox{\boldmath $\xi$}\in{\cal N}(x,\sigma,\mbox{\boldmath $a$})\,\right\}.
$$
This is analogous to the reconstruction formula  (2.9) of  $\sigma(]0,\,1])$ for $\sigma\in N_x$ via the 
set of numerical values:
$$\displaystyle
\left\{\Vert\text{Sym}\,\nabla\mbox{\boldmath $v$}_n\Vert_{L^2(B\cap\Omega)}\,\vert\,
\mbox{\boldmath $\xi$}=\{(\mbox{\boldmath $v$}_n,q_n)\}\in{\cal N}(x,\sigma,\mbox{\boldmath $a$}),
B\in\cup_{z\in\Omega} \beta_z\,\right\}.
$$

\subsection{The singular sources method via the probe method}

 First, considering the meaning of the original singular sources method as described in \cite{IPS},
 we introduce the indicator function for the singular sources method for the Stokes system.
 
 $\quad$
 
 {\bf\noindent Definition 3.3.}  Fix  $\mbox{\boldmath $a$}\not=\mbox{\boldmath $0$}$.
 We call the function defined by
 $$\displaystyle
 \Omega\setminus\overline{D}\ni x\longmapsto \mbox{\boldmath $w$}_x(x)\mbox{\boldmath $a$}\cdot\mbox{\boldmath $a$},
 $$ 
 the  indicator function of the singular sources method.

 \subsubsection{The third indicator function or IPS function}

{\bf\noindent Definition 3.4.}
Let the pair $(\mbox{\boldmath $w$},q)\in H^2(\Omega\setminus\overline{D})\times H^1(\Omega\setminus\overline{D})$ be
the solution of
$$\left\{
 \begin{array}{ll}
 \mu\Delta\mbox{\boldmath $w$}-\nabla q=\mbox{\boldmath $0$}, & z\in\Omega\setminus\overline{D},
 \\
 \\
 \displaystyle
 \nabla\cdot\mbox{\boldmath $w$}=0, &  z\in\Omega\setminus\overline{D},\\
 \\
 \displaystyle
 \mbox{\boldmath $w$}=\mbox{\boldmath $0$}, & z\in\partial D,
 \\
 \\
 \displaystyle
 \mbox{\boldmath $w$}=\mbox{\boldmath $J$}(z-x)
 \mbox{\boldmath $a$}, & z\in\partial\Omega.
 \end{array}
 \right.
 \tag {3.10}
 $$
 As well as $\mbox{\boldmath $w$}_x$ there exists the unique
 $M_3(\Bbb R)$-valued function $\mbox{\boldmath $w$}_x^1$ such that
 $$\displaystyle
 \mbox{\boldmath $w$}=\mbox{\boldmath $w$}_x^1\mbox{\boldmath $a$}.
 $$
 The $q$ is unique modulo a constant function.  We denote it by $q_x^1=q_x^1(\,\cdot\,;\mbox{\boldmath $a$})$.

 $\quad$

{\bf\noindent Definition 3.5.}  
 Given $x\in\Omega\setminus\overline{D}$, define
$$\begin{array}{ll}
\displaystyle
\mbox{\boldmath $W$}_x(y)=\mbox{\boldmath $w$}_x(y)+\mbox{\boldmath $w$}_x^1(y),
&
\displaystyle
y\in\Omega\setminus\overline{D}.
\end{array}
$$
We call this $M_3(\Bbb R^3)$-valued function the lifting of  {\it third indicator function} or IPS function 
defined below.

$\quad$

Note that, for each constant vector $\mbox{\boldmath $a$}$
the function $\mbox{\boldmath $W$}(z)=\mbox{\boldmath $W$}_x(z)\mbox{\boldmath $a$}$ satisfies
$$\left\{
 \begin{array}{ll}
 \mu\Delta\mbox{\boldmath $W$}-\nabla q=\mbox{\boldmath $0$}, & z\in\Omega\setminus\overline{D},
 \\
 \\
 \displaystyle
 \nabla\cdot\mbox{\boldmath $W$}=0, &  z\in\Omega\setminus\overline{D},\\
 \\
 \displaystyle
 \mbox{\boldmath $W$}=-\mbox{\boldmath $J$}(z-x)\mbox{\boldmath $a$}, & z\in\partial D,
 \\
 \\
 \displaystyle
 \mbox{\boldmath $W$}=\mbox{\boldmath $J$}(z-x)
 \mbox{\boldmath $a$}, & z\in\partial\Omega,
 \end{array}
 \right.
 $$
where $q=q_x+q_x^1$ (see Definitions 3.1 and 3.4).
This is the system which characterizes the lifting of the IPS function.

$\quad$

 {\bf\noindent Definition 3.6.}  Fix $\mbox{\boldmath $a$}\not=\mbox{\boldmath $0$}$.  We call the function defined by
 $$\displaystyle
 \Omega\setminus\overline{D}\ni x\longmapsto \mbox{\boldmath $W$}_x(x)
 \mbox{\boldmath $a$}\cdot\mbox{\boldmath $a$},
 $$
 the {\it third indicator function} or {\it the IPS function}.

 $\quad$

We have the trivial decomposition
$$\displaystyle
\mbox{\boldmath $W$}_x(x)\mbox{\boldmath $a$}\cdot\mbox{\boldmath $a$}
=\mbox{\boldmath $w$}_x(x)\mbox{\boldmath $a$}\cdot\mbox{\boldmath $a$}
+\mbox{\boldmath $w$}_x^1(x)\mbox{\boldmath $a$}\cdot\mbox{\boldmath $a$}.
\tag {3.11}
$$
This is a point-wise decomposition.
We show that there is another decomposition of  the IPS function.

For the purpose we introduce 

$\quad$

{\bf\noindent Definition 3.7.}
Let $\mbox{\boldmath $a$}$ and $\mbox{\boldmath $b$}$ be arbitrary vectors.
Let $(x,y)\in (\Omega\setminus\overline{D})^2$.  
Define
$$\begin{array}{l}
\displaystyle
\,\,\,\,\,\,
I^1(x,\mbox{\boldmath $a$};y,\mbox{\boldmath $b$})\\
\\
\displaystyle
=\int_{\Omega\setminus\overline{D}}2\mu\,\text{Sym}\nabla(\mbox{\boldmath $w$}_x^1(z)\mbox{\boldmath $a$})
\cdot\text{Sym}\nabla(\mbox{\boldmath $w$}_y^1(z)\mbox{\boldmath $b$})\,dz
\\
\\
\displaystyle
\,\,\,
+
\int_{\Bbb R^3\setminus\overline{\Omega}}
2\mu\,\text{Sym}\,\nabla(\mbox{\boldmath $J$}(z-x)\mbox{\boldmath $a$})
\cdot\text{Sym}\,\nabla(\mbox{\boldmath $J$}(z-y)\mbox{\boldmath $b$})\,dz.
\end{array}
$$
and 
$$\displaystyle
I^1(x,\mbox{\boldmath $a$})\equiv I^1(x,\mbox{\boldmath $a$};x,\mbox{\boldmath $a$}).
$$

$\quad$

We have also the trivial symmetry relation
$$\displaystyle
I^1(x,\mbox{\boldmath $a$};y,\mbox{\boldmath $b$})=I^1(y,\mbox{\boldmath $b$};x,\mbox{\boldmath $a$}).
$$
 Besides, it follows from (1.2) with $\mbox{\boldmath $f$}=\mbox{\boldmath $J$}(\,\cdot\,-x)\mbox{\boldmath $a$}$
 and $\mbox{\boldmath $\phi$}=\mbox{\boldmath $w$}_y^1\mbox{\boldmath $b$}$, we have
 $$\begin{array}{ll}
 \displaystyle
 <\Lambda_D(\mbox{\boldmath $J$}(\,\cdot\,-x)\mbox{\boldmath $a$}\vert_{\partial\Omega}),\mbox{\boldmath $J$}(\,\cdot\,-y)\mbox{\boldmath $b$}\vert_{\partial\Omega}>
 &
 \displaystyle
 =\int_{\Omega\setminus\overline{D}}2\mu\,\text{Sym}\nabla(\mbox{\boldmath $w$}_x^1(z)\mbox{\boldmath $a$})
\cdot\text{Sym}\nabla(\mbox{\boldmath $w$}_y^1(z)\mbox{\boldmath $b$})\,dz.
 \end{array}
 $$
 Thus we obtain the expression
 $$\begin{array}{ll}
 \displaystyle
 I^1(x,\mbox{\boldmath $a$};y,\mbox{\boldmath $b$})
 &
 \displaystyle
 =<\Lambda_D(\mbox{\boldmath $J$}(\,\cdot\,-x)\mbox{\boldmath $a$}\vert_{\partial\Omega}),
 \mbox{\boldmath $J$}(\,\cdot\,-y)\mbox{\boldmath $b$}\vert_{\partial\Omega}>
 \\
 \\
 \displaystyle
 &
 \displaystyle
 \,\,\,
 +\int_{\Bbb R^3\setminus\overline{\Omega}}
2\mu\,\text{Sym}\,\nabla(\mbox{\boldmath $J$}(z-x)\mbox{\boldmath $a$})
\cdot\text{Sym}\,\nabla(\mbox{\boldmath $J$}(z-y)\mbox{\boldmath $b$})\,dz.
\end{array}
\tag {3.12}
 $$
This means that the $I^1(x,\mbox{\boldmath $a$};y,\mbox{\boldmath $b$})$ can be calculated from $\Lambda_D$
without making use of needle sequences.

The result mentioned below is the core of the integrated theory of the probe and singular sources methods.

 \proclaim{\noindent Theorem 3.3.}
 Let $\mbox{\boldmath $a$}$ and $\mbox{\boldmath $b$}$ be arbitrary vectors.
 We have
$$\begin{array}{l}
\displaystyle
\,\,\,\,\,\,
\mbox{\boldmath $W$}_x(y)\mbox{\boldmath $a$}\cdot\mbox{\boldmath $b$}
\\
\\
\displaystyle
=\int_{\partial\Omega}\sigma(\mbox{\boldmath $w$}_x(z)\mbox{\boldmath $a$},q_x(z;\mbox{\boldmath $a$}))\mbox{\boldmath $\nu$}\cdot
\mbox{\boldmath $J$}(z-y)\mbox{\boldmath $b$}\,dz
-\int_{\partial\Omega}\mbox{\boldmath $J$}(z-x)\mbox{\boldmath $a$}
\cdot\sigma(\mbox{\boldmath $w$}_y(z)\mbox{\boldmath $b$},q_y(z;\mbox{\boldmath $b$}))
\mbox{\boldmath $\nu$}\,dS(z)
\\
\\
\displaystyle
\,\,\,
+I(x,\mbox{\boldmath $a$};y,\mbox{\boldmath $b$})+I^1(x,\mbox{\boldmath $a$};y,\mbox{\boldmath $b$}).
\end{array}
\tag {3.13}
$$
In particular, we have
$$\displaystyle
\mbox{\boldmath $W$}_x(x)\mbox{\boldmath $a$}\cdot\mbox{\boldmath $a$}
=I(x,\mbox{\boldmath $a$})+I^1(x,\mbox{\boldmath $a$}).
\tag {3.14}
$$
\endproclaim

It follows from (3.13) that
$$\displaystyle
\frac{1}{2}\left(\mbox{\boldmath $W$}_x(y)+\mbox{\boldmath $W$}_y(x)^T\,\right)\,\mbox{\boldmath $a$}\cdot\mbox{\boldmath $b$}
=I(x,\mbox{\boldmath $a$};y,\mbox{\boldmath $b$})+I^1(x,\mbox{\boldmath $a$};y,\mbox{\boldmath $b$}).
\tag {3.15}
$$
The left hand-side on this equation is the symmetrization of  the function
$\mbox{\boldmath $W$}_x(y)\mbox{\boldmath $a$}\cdot\mbox{\boldmath $b$}$ with respect to
two sets of variables $(x,\mbox{\boldmath $a$})$ and $(y,\mbox{\boldmath $b$})$.

Therefore from Proposition 3.2, (3.12) and (3.15), one can calculate the symmetric part of  the lifting of the third indicator function by using $\Lambda_0-\Lambda_D$ acting needle sequences and $\Lambda_D$.

$\quad$

{\bf\noindent Remark 3.2.}
Note also that from (3.3) of Proposition 3.1 and (3.13) we have the expression
$$\begin{array}{ll}
\displaystyle
\frac{1}{2}\left(\mbox{\boldmath $W$}_x(y)-\mbox{\boldmath $W$}_y(x)^T\,\right)\,\mbox{\boldmath $a$}\cdot\mbox{\boldmath $b$}
&
\displaystyle
=\lim_{n\rightarrow\infty}
<(\Lambda_0-\Lambda_D)\mbox{\boldmath $v$}_n\vert_{\partial\Omega},
\mbox{\boldmath $J$}(\,\cdot\,-y)\mbox{\boldmath $b$}\vert_{\partial\Omega}>
\\
\\
\displaystyle
&
\displaystyle
\,\,\,
-\lim_{n\rightarrow\infty}
<(\Lambda_0-\Lambda_D)\mbox{\boldmath $v$}_n'\vert_{\partial\Omega},
\mbox{\boldmath $J$}(\,\cdot\,-x)\mbox{\boldmath $a$}\vert_{\partial\Omega}>,
\end{array}
$$
where $\{(\mbox{\boldmath $v$}_n, q_n)\}\in{\cal N}(x,\sigma,\mbox{\boldmath $a$})$,
 $\{(\mbox{\boldmath $v$}_n', q_n')\}\in{\cal N}(y,\sigma',\mbox{\boldmath $b$})$, $\sigma\in N_x$ and $\sigma'\in N_y$
 and satisfy $\sigma\cap\overline{D}=\sigma'\cap\overline{D'}=\emptyset$.
Therefore, the value of $\mbox{\boldmath $W$}_x(y)\mbox{\boldmath $a$}\cdot\mbox{\boldmath $b$}$ 
at $(x,y)\in(\Omega\setminus\overline{D})^2$ itself
has an expression in terms of  $\Lambda_0-\Lambda_D$, $\Lambda_D$ together with needle sequences.

$\quad$

Next we are interested in the behaviour of the third indicator function $\mbox{\boldmath $W$}_x(x)\mbox{\boldmath $a$}\cdot\mbox{\boldmath $a$}$ together with functions 
$\mbox{\boldmath $w$}_x(x)
\mbox{\boldmath $a$}\cdot\mbox{\boldmath $a$}$ and $\mbox{\boldmath $w$}_x^1(x)
\mbox{\boldmath $a$}\cdot\mbox{\boldmath $a$}$ as $x$ approaches points on $\partial\Omega$ and $\partial D$.
By the well-posedness of (3.1), (3.10) and the Sobolev imbedding we have, for each fixed positive number $\epsilon$
 $$\displaystyle
 \sup_{x\in\Omega\setminus\overline{D},\,\text{dist}\,(x,\partial D)>\epsilon}\,\vert \mbox{\boldmath $w$}_x(x)\mbox{\boldmath $a$}
 \cdot\mbox{\boldmath $a$}\vert<\infty
 \tag {3.16}
 $$
 and
 $$\displaystyle
 \sup_{x\in\Omega\setminus\overline{D},\,\text{dist}\,(x,\partial\Omega)>\epsilon}\,\vert \mbox{\boldmath $w$}_x^1(x)
 \mbox{\boldmath $a$}\cdot\mbox{\boldmath $a$}\vert<\infty.
 \tag {3.17}
 $$
 And also we have
 $$\displaystyle
 0\le\sup_{x\in\Omega\setminus\overline{D},\,\text{dist}\,(x,\partial\Omega)>\epsilon}
 \,I^1(x,\mbox{\boldmath $a$})<\infty.
 \tag {3.18}
 $$

 Now we are ready to state the behaviour of the third indicator function on the outer surface $\partial\Omega$ and 
 obstacle surface $\partial D$.

 \proclaim{\noindent Theorem 3.4.}
 Let $\mbox{\boldmath $a$}\not=\mbox{\boldmath $0$}$.
 We have

 \noindent
 (i)  given an arbitrary point $c\in\partial\Omega$ 
 $$\displaystyle
 \lim_{x\rightarrow c}\mbox{\boldmath $W$}_x(x)\mbox{\boldmath $a$}\cdot\mbox{\boldmath $a$}=\infty.
 $$

 \noindent
 (ii)   given an arbitrary point $d\in\partial D$

 \noindent
 $$\displaystyle
 \lim_{x\rightarrow d}\mbox{\boldmath $W$}_x(x)\mbox{\boldmath $a$}\cdot\mbox{\boldmath $a$}=-\infty.
 $$

 \noindent
 (iii)  Let $\epsilon_i>0$, $i=1,2$.  Then we have
 $$\displaystyle
 \sup_{\epsilon_1<\text{dist}\,(x,\partial\Omega),\,\epsilon_2<\text{dist}\,(x,\partial D)\,}
 \vert\mbox{\boldmath $W$}_x(x)\mbox{\boldmath $a$}\cdot\mbox{\boldmath $a$}\vert<\infty.
 $$

 \endproclaim

 {\it\noindent Proof.}  
 The proof of (i) is as follows.  By (3.7) of (iii) in Theorem 3.1 and (3.14)
 we have, for all $x\in\Omega\setminus\overline{D}$ with $\text{dist}\,(x,\partial D)>\epsilon$
 $$\displaystyle
 \mbox{\boldmath $W$}_x(x)\mbox{\boldmath $a$}\cdot\mbox{\boldmath $a$}
 >I^1(x,\mbox{\boldmath $a$})-C,
 $$
 where $0<\epsilon<<1$ and $C>>1$ and they are independent of $x$.
 By the expression (3.12) we have
 $$\displaystyle
 I^1(x,\mbox{\boldmath $a$})
\ge \int_{\Bbb R^3\setminus\overline{\Omega}}
2\mu\vert\text{Sym}\,\nabla(\mbox{\boldmath $J$}(z-x)\mbox{\boldmath $a$})\vert^2\,dz.
$$
Applying Proposition 2.2 with a suitable finite cone $V$ with vertex at $c\in\partial\Omega$
with the help of Fatou's lemma, we conclude that the integral on this right-hand side blows up
as $x\rightarrow c$.
 
 \noindent
 (ii) is a direct consequence of  the expression  (3.14), estimate (3.18)
 and (3.6) of (ii) in Theorem 3.1.
 
 \noindent
 (iii) is a consequence of  the trivial expression (3.11), estimates (3.16) and (3.17).
 
 \noindent
 $\Box$

 And rewrite (3.11) as
 $$\displaystyle
 \mbox{\boldmath $w$}_x(x)\mbox{\boldmath $a$}\cdot\mbox{\boldmath $a$}
 = \mbox{\boldmath $W$}_x(x)\mbox{\boldmath $a$}\cdot\mbox{\boldmath $a$}
 -\mbox{\boldmath $w$}_x^1(x)\mbox{\boldmath $a$}\cdot\mbox{\boldmath $a$}.
 $$
 From this together with estimate (3.17) and (ii) of Theorem 3.4 
 we obtain the blowing up property
 of the obstacle surface of the indicator function for the singular sources method.

\proclaim
 {\bf\noindent Corollary 3.1.}
 Let $\mbox{\boldmath $a$}\not=\mbox{\boldmath $0$}$.
 Given an arbitrary point $d\in\partial D$ we have
 $$\displaystyle
 \lim_{x\rightarrow d}\mbox{\boldmath $w$}_x(x)\mbox{\boldmath $a$}\cdot\mbox{\boldmath $a$}=-\infty.
 \tag {3.19}
 $$

 \endproclaim
 
\noindent
The point of our approach is: in the proof of the blowing up property (3.19), we have nothing to make use of the integral equation approach in the previous Potthast's singular sources method.
Instead we use only the consequence of the decompositions (3.11) and (3.14), that is, the equation
$$\displaystyle
\mbox{\boldmath $w$}_x(x)\mbox{\boldmath $a$}\cdot\mbox{\boldmath $a$}
+\mbox{\boldmath $w$}_x^1(x)\mbox{\boldmath $a$}\cdot\mbox{\boldmath $a$}
=I(x,\mbox{\boldmath $a$})+I^1(x,\mbox{\boldmath $a$}).
$$

The next problem is to give a calculation formula of indicator function
$\mbox{\boldmath $w$}_x(x)\mbox{\boldmath $a$}\cdot\mbox{\boldmath $a$}$.
We follow the idea in \cite{IPS}.

$\quad$

 {\bf\noindent Definition 3.8.}  Fix $\mbox{\boldmath $a$}\not=\mbox{\boldmath $0$}$.
 Let $x\in\Omega$ and $\sigma\in N_x$.  Given a needle sequence $\mbox{\boldmath $\xi$}=\{(\mbox{\boldmath $v$}_n, q_n)\}\in{\cal N}(x,\sigma,\mbox{\boldmath $a$})$
define
$$\left\{
\begin{array}{ll}
\displaystyle
\mbox{\boldmath $J$}_n(z;\mbox{\boldmath $\xi$})=\mbox{\boldmath $J$}(z-x)\mbox{\boldmath $a$}-\mbox{\boldmath $v$}_n(z), & z\in\Omega\setminus\{x\},
\\
\\
\displaystyle
p_n(z;\mbox{\boldmath $\xi$})=\mbox{\boldmath $p$}(z-x)\cdot\mbox{\boldmath $a$}
-q_n(z), & z\in\Omega\setminus\{x\}.
\end{array}
\right.
\tag {3.20}
$$

 $\quad$

 Here we describe two important properties of the pair  $(\mbox{\boldmath $J$}_n(\,\cdot\,;\mbox{\boldmath $\xi$}), p_n(\,\cdot\,;\mbox{\boldmath $\xi$}))$
 which is coming from Definition 3.1.
 
 $\bullet$  The pair $(\mbox{\boldmath $J$}_n(\,\cdot\,;\mbox{\boldmath $\xi$}), p_n(\,\cdot\,;\mbox{\boldmath $\xi$}))$ satisfies the Stokes system in $\Omega\setminus\{x\}$ and their singularity at $z=x$ coincides with
 that of $(\mbox{\boldmath $J$}(\,\cdot-x)\mbox{\boldmath $a$}, \mbox{\boldmath $p$}(\,\cdot\,-x)\cdot\mbox{\boldmath $a$})$.
 
 $\bullet$  If $\sigma\cap\overline{D}=\emptyset$, then the sequence
 $\{(\mbox{\boldmath $J$}_n(\,\cdot\,;\mbox{\boldmath $\xi$}), p_n(\,\cdot\,;\mbox{\boldmath $\xi$}))\}$ 
 converges to $\{(\mbox{\boldmath $0$},0)\}$ in $H^2(D)\times H^1(D)$.

 These properties combined with integration by parts yield the calculation formula
 of indicator function $\mbox{\boldmath $w$}_x(x)\mbox{\boldmath $a$}\cdot\mbox{\boldmath $a$}$.

\proclaim{\noindent Theorem 3.5.}
Fix $\mbox{\boldmath $a$}\not=\mbox{\boldmath $0$}$.
Let $x\in\Omega\setminus\overline{D}$ and $\sigma\in N_x$.
If $\sigma\cap\overline{D}=\emptyset$, then, for all needle sequences 
$\mbox{\boldmath $\xi$}=\{(\mbox{\boldmath $v$}_n, q_n)\}\in{\cal N}(x,\sigma,\mbox{\boldmath $a$})$ we have
$$\displaystyle
-\lim_{n\rightarrow\infty}<(\Lambda_0-\Lambda_D)\mbox{\boldmath $v$}_n\vert_{\partial\Omega},
\mbox{\boldmath $J$}_n(\,\cdot\,;\mbox{\boldmath $\xi$})\vert_{\partial\Omega}>
=\mbox{\boldmath $w$}_x(x)\mbox{\boldmath $a$}\cdot\mbox{\boldmath $a$}.
\tag {3.21}
$$

\endproclaim

Thus we have established the Side A of the singular sources method which is based on (3.16), Corollary 3.1 and Theorem 3.5.  The proof is given in Subsection 4.3.

However, as described in \cite{IPS}, the Side B for this singular sources method was not known.    It means that
 when $x\in\overline{D}$ or $x\in\Omega\setminus\overline{D}$ and $\sigma\in N_x$ satisfies $\sigma\cap\overline{D}\not=\emptyset$ the behaviour of the sequence
 $$\displaystyle
 -<(\Lambda_0-\Lambda_D)\mbox{\boldmath $v$}_n\vert_{\partial\Omega},
\mbox{\boldmath $J$}_n(\,\cdot\,;\mbox{\boldmath $\xi$})\vert_{\partial\Omega}>,
$$
 was not clear.  In the next subsection we reformulate the probe method
 as the completely integrated method in the sense that the indicator function for  the reformulated method
 has the analogous expression to that of the singular sources method.
 As a byproduct we have succeeded to show that the singular sources method
 introduced in this section also has the Side B.  It means that the asymptotic behaviour
 of the  sequence above is clarified when $x\in\overline{D}$ or $x\in\Omega\setminus\overline{D}$ and
 $\sigma\in N_x$ satisfies $\sigma\cap\overline{D}\not=\emptyset$.  This is a new fact
 not mentioned in \cite{IPS} even for the Laplace equation case.

\subsection{Completely integrated method}

In this section we present the result  that explains the origin of  the naming 
`the integrated theory of  the probe and singular sources methods'.

The idea is just to replace the singular solution of the probe method with another one which keeps the original singularity.
Needless to say, such a replacement is not unique.  Here we choose the singular solution in such a way that
its trace on the surface $\partial\Omega$ vanishes.  That is the Green function.
As pointed out in the previous works on the probe method, the important things in choosing the singular solution
is only the strength of the singularity.  So based on another singular solution, one can define the needle sequence and construct the theory of the probe method which has both the Side A and Side B.  Then corresponding to the new singular solution, one can define
the reflected solution and construct the new indicator function for the singular sources method in terms of the probe method just mentioned above.  In this time, the third indicator function coincides with all the indicator functions.
So we obtain the singular sources method having not only the Side A but also Side B.

Now let us explain more precisely step by step.

\subsubsection{Side A}

First we introduce another singular solution of the Stokes system for a technical reason.

$\quad$

{\bf\noindent Definition 3.9.}  Fix $\mbox{\boldmath $a$}\not=\mbox{\boldmath $0$}$.
Given $x\in\Omega$, let $(\mbox{\boldmath $\epsilon$}, \rho)
=(\mbox{\boldmath $\epsilon$}_x(\,\cdot\,;\mbox{\boldmath $a$}), \rho_x(\,\cdot\,;\mbox{\boldmath $a$}))\in H^2(\Omega)\times H^1(\Omega)$ be the solution of
 $$\left\{
 \begin{array}{ll}
 \mu\Delta\mbox{\boldmath $\epsilon$}-\nabla \rho=\mbox{\boldmath $0$}, & z\in\Omega,
 \\
 \\
 \displaystyle
 \nabla\cdot\mbox{\boldmath $\epsilon$}=0, &  z\in\Omega,
 \\
 \\
 \displaystyle
 \mbox{\boldmath $\epsilon$}=-\mbox{\boldmath $J$}(z-x)\mbox{\boldmath $a$}, & z\in\partial\Omega.
 \end{array}
 \right.
 \tag {3.22}
 $$
Define
$$
\begin{array}{ll}
\displaystyle
\mbox{\boldmath $J$}_{\Omega}(z;x,\mbox{\boldmath $a$})=\mbox{\boldmath $J$}(z-x)\mbox{\boldmath $a$}+\mbox{\boldmath $\epsilon$}_x(z;\mbox{\boldmath $a$}),
& z\in\Omega\setminus\{x\}
\end{array}
\tag {3.23}
$$
and
$$
\begin{array}{ll}
\displaystyle
p_{\Omega}(z;x,\mbox{\boldmath $a$})=\mbox{\boldmath $p$}(z-x)\cdot\mbox{\boldmath $a$}+\rho_x(z;\mbox{\boldmath $a$}),
& z\in\Omega\setminus\{x\}.
\end{array}
\tag {3.24}
$$

$\quad$

 By the well-posedness of (3.22), we have the estimate
 $$\displaystyle
 \forall\epsilon>0 \,\sup_{x\in\Omega,\,\text{dist}\,(x,\partial\Omega)>\epsilon}
 \Vert\mbox{\boldmath $\epsilon$}_x\Vert_{H^2(\Omega\setminus\overline{D})}<\infty.
 \tag {3.25}
 $$

Note that $\mbox{\boldmath $J$}_{\Omega}(\,\cdot\,;x,\mbox{\boldmath $a$})=\mbox{\boldmath $0$}$ on $\partial\Omega$.
Besides the pair $(\mbox{\boldmath $J$}_{\Omega},p_{\Omega})=(\mbox{\boldmath $J$}_{\Omega}(\,\cdot\,;x,\mbox{\boldmath $a$}), p_{\Omega}(\,\cdot\,;x,\mbox{\boldmath $a$}))$
satisfies, in the sense of distributions
$$\left\{
 \begin{array}{ll}
 \mu\Delta\mbox{\boldmath $J$}_{\Omega}-\nabla p+\mbox{\boldmath $a$}
 \delta(z-x)=\mbox{\boldmath $0$}, & z\in\Omega,
 \\
 \\
 \displaystyle
 \nabla\cdot\mbox{\boldmath $J$}_{\Omega}=0, &  z\in\Omega.
 \end{array}
 \right.
 $$
The leading singularity of this pair completely coincides with the Stokeslet.

Based on the pair $(\mbox{\boldmath $J$}_{\Omega},p_{\Omega})$ we introduce another reflected solution for the Stokes system.

$\quad$

 {\bf\noindent Definition 3.10.}
 Given $x\in\Omega\setminus\overline{D}$, let the pair $(\mbox{\boldmath $w$}, q)=(\mbox{\boldmath $w$}_x^*(\,\cdot\,;\mbox{\boldmath $a$}), q_x^*(\,\cdot\,;\mbox{\boldmath $a$}))$ be the solution of
  $$\left\{
 \begin{array}{ll}
 \mu\Delta\mbox{\boldmath $w$}-\nabla q=\mbox{\boldmath $0$}, & z\in\Omega\setminus\overline{D},
 \\
 \\
 \displaystyle
 \nabla\cdot\mbox{\boldmath $w$}=0, &  z\in\Omega\setminus\overline{D},
 \\
 \\
 \displaystyle
 \mbox{\boldmath $w$}=-\mbox{\boldmath $J$}_{\Omega}(z;x,\mbox{\boldmath $a$}), & z\in\partial D,
 \\
 \\
 \displaystyle
 \mbox{\boldmath $w$}=\mbox{\boldmath $0$}, & z\in\partial\Omega.
 \end{array}
 \right.
 $$

  $\quad$

  By Definition 3.1, (3.23) and Definition 3.9, we have the expression
  $$\displaystyle
  \mbox{\boldmath $w$}_x^*(z;\mbox{\boldmath $a$})=\mbox{\boldmath $w$}_x(z)\mbox{\boldmath $a$}+
  \mbox{\boldmath $z$}_x(z;\mbox{\boldmath $a$}),
  \tag {3.26}
  $$
  where $\mbox{\boldmath $z$}=\mbox{\boldmath $z$}_x(z;\mbox{\boldmath $a$})$ with a $r=r_x(\,\cdot\,;\mbox{\boldmath $a$})$ solves
  $$\left\{
 \begin{array}{ll}
 \mu\Delta\mbox{\boldmath $z$}-\nabla r=\mbox{\boldmath $0$}, & z\in\Omega\setminus\overline{D},
 \\
 \\
 \displaystyle
 \nabla\cdot\mbox{\boldmath $z$}=0, &  z\in\Omega\setminus\overline{D},
 \\
 \\
 \displaystyle
 \mbox{\boldmath $z$}=-\mbox{\boldmath $\epsilon$}_x(z;\mbox{\boldmath $a$}), & z\in\partial D,
 \\
 \\
 \displaystyle
 \mbox{\boldmath $z$}=\mbox{\boldmath $0$}, & z\in\partial\Omega.
 \end{array}
 \right.
 $$
 By the well-posedness, we have the estimate
 $$\displaystyle
 \forall\epsilon>0 \,\sup_{x\in\Omega,\,\text{dist}(x,\partial\Omega)>\epsilon}
 \Vert\mbox{\boldmath $z$}_x\Vert_{H^2(\Omega\setminus\overline{D})}<\infty.
 \tag {3.27}
  $$

  $\quad$

  {\bf\noindent Definition 3.11.}  Fix $\mbox{\boldmath $a$}\not=\mbox{\boldmath $0$}$.
  Define
  $$\begin{array}{ll}
  \displaystyle
  I^*(x,\mbox{\boldmath $a$})
  =-\Vert 2\mu\,\text{Sym}\,\nabla\mbox{\boldmath $w$}_x^*(\,\cdot\,;\mbox{\boldmath $a$})\Vert_{L^2(\Omega\setminus\overline{D})}^2
  -\Vert 2\mu\,\text{Sym}\,\nabla \mbox{\boldmath $J$}_{\Omega}(\,\cdot\,;x,\mbox{\boldmath $a$})\Vert_{L^2(D)}^2,
  &
  \displaystyle
  x\in\Omega\setminus\overline{D}.
  \end{array}
  $$

  $\quad$

Now we describe a result which shows that the function $I^*(x,\mbox{\boldmath $a$})$ and the sequence
$\{<(\Lambda_0-\Lambda_D)\mbox{\boldmath $J$}_n(\,\cdot\,;\mbox{\boldmath $\xi$})\vert_{\partial\Omega},
\mbox{\boldmath $J$}_n(\,\cdot\,;\mbox{\boldmath $\xi$})\vert_{\partial\Omega}>\}$
for $\mbox{\boldmath $\xi$}\in{\cal N}(x,\sigma,\mbox{\boldmath $a$})$
plays almost the same role of the indicator function and sequence for the original probe method in the Side A.

\proclaim{\noindent Theorem 3.6.}  Fix $\mbox{\boldmath $a$}\not=\mbox{\boldmath $0$}$.

\noindent
(i)  Let $x\in\Omega\setminus\overline{D}$ and $\sigma\in N_x$.
If $\sigma\cap\overline{D}=\emptyset$, then,
for all needle sequences $\mbox{\boldmath $\xi$}=\{(\mbox{\boldmath $v$}_n,q_n)\}\in{\cal N}(x,\sigma,\mbox{\boldmath $a$})$ we have
$$\displaystyle
\lim_{n\rightarrow\infty}
<(\Lambda_0-\Lambda_D)\mbox{\boldmath $J$}_n(\,\cdot\,;\mbox{\boldmath $\xi$})\vert_{\partial\Omega},
\mbox{\boldmath $J$}_n(\,\cdot\,;\mbox{\boldmath $\xi$})\vert_{\partial\Omega}>
=I^*(x,\mbox{\boldmath $a$}).
\tag {3.28}
$$

\noindent
(ii)  We have
$$\displaystyle
\lim_{x\rightarrow a\in\partial D}I^*(x,\mbox{\boldmath $a$})=-\infty.
\tag {3.29}
$$

\noindent
(iii)  For each positive numbers $\epsilon_1$ and $\epsilon_2$ it holds that
$$\displaystyle
-\infty<\inf_{\text{dist}\,(x,\partial D)>\epsilon_1, \text{dist}\,(x,\partial\Omega)>\epsilon_2}\,I^*(x,\mbox{\boldmath $a$})<0.
\tag {3.30}
$$

 \endproclaim
{\it\noindent Proof.}  First we give a proof of (i).
For $z\in\partial\Omega$, we have
$$\begin{array}{ll}
\displaystyle
\mbox{\boldmath $J$}_n(z;\mbox{\boldmath $\xi$})
&
\displaystyle
=\mbox{\boldmath $J$}(z-x)\mbox{\boldmath $a$}-\mbox{\boldmath $v$}_n(z)
\\
\\
\displaystyle
&
\displaystyle
=-(\mbox{\boldmath $\epsilon$}_x(z;\mbox{\boldmath $a$})+\mbox{\boldmath $v$}_n(z))
\end{array}
$$
and thus
$$\begin{array}{l}
\displaystyle
\,\,\,\,\,\,
<(\Lambda_0-\Lambda_D)\mbox{\boldmath $J$}_n(\,\cdot\,;\mbox{\boldmath $\xi$})\vert_{\partial\Omega},
\mbox{\boldmath $J$}_n(\,\cdot\,;\mbox{\boldmath $\xi$})\vert_{\partial\Omega}>
\\
\\
\displaystyle
=<(\Lambda_0-\Lambda_D)
(\mbox{\boldmath $\epsilon$}_x(\,\cdot\,;\mbox{\boldmath $a$})+\mbox{\boldmath $v$}_n)
\vert_{\partial\Omega},
(\mbox{\boldmath $\epsilon$}_x(\,\cdot\,;\mbox{\boldmath $a$})+\mbox{\boldmath $v$}_n)
\vert_{\partial\Omega}>.
\end{array}
\tag {3.31}
$$
Here one ses that the pair 
$(\mbox{\boldmath $\epsilon$}_x(\,\cdot\,;\mbox{\boldmath $a$})+\mbox{\boldmath $v$}_n,\rho_x(\,\cdot\,;\mbox{\boldmath $a$})+q_n)$
becomes the needle sequence for $(x,\sigma)$ based on the pair of the singular solutions $(\mbox{\boldmath $J$}(\,\cdot\,-x)\mbox{\boldmath $a$},\mbox{\boldmath $p$}(\,\cdot\,-x)\cdot\mbox{\boldmath $a$})$ in Definition 2.1 replaced with
$(\mbox{\boldmath $J$}_{\Omega}(\,\cdot\,;x,\mbox{\boldmath $a$}),p_{\Omega}(\,\cdot\,;x,\mbox{\boldmath $a$}))$ given by (3.23) and (3.24).
Thus by (3.3) we have
$$\begin{array}{l}
\displaystyle
\,\,\,\,\,\,
<(\Lambda_0-\Lambda_D)
(\mbox{\boldmath $\epsilon$}_x(\,\cdot\,;\mbox{\boldmath $a$})+\mbox{\boldmath $v$}_n)
\vert_{\partial\Omega},
(\mbox{\boldmath $\epsilon$}_x(\,\cdot\,;\mbox{\boldmath $a$})+\mbox{\boldmath $v$}_n)
\vert_{\partial\Omega}>
\\
\\
\displaystyle
=-\Vert 2\mu\,\text{Sym}\,\nabla\mbox{\boldmath $w$}_n\Vert_{L^2(\Omega\setminus\overline{D})}^2
-\Vert 2\mu\,\text{Sym}\,\nabla(\mbox{\boldmath $\epsilon$}_x(\,\cdot
\,;\mbox{\boldmath $a$})+\mbox{\boldmath $v$}_n)
\Vert_{L^2(D)}^2,
\end{array}
\tag {3.32}
$$
where the $\mbox{\boldmath $w$}_n=\mbox{\boldmath $w$}$ with a $q_n=q$ satisfies the Stokes system (3.1) with 
$-\mbox{\boldmath $J$}(z-x)\mbox{\boldmath $a$}$ for $z\in\partial D$ replaced with $-(\mbox{\boldmath $\epsilon$}_x(z;\mbox{\boldmath $a$})+\mbox{\boldmath $v$}_n(z;\sigma,\mbox{\boldmath $a$}))$.
Since we have in $H^{2}(D)$
$$\displaystyle
\mbox{\boldmath $\epsilon$}_x(z;\mbox{\boldmath $a$})+\mbox{\boldmath $v$}_n(z;\sigma,\mbox{\boldmath $a$})
\rightarrow
\mbox{\boldmath $\epsilon$}_x(z;\mbox{\boldmath $a$})+\mbox{\boldmath $J$}(z-x)\mbox{\boldmath $a$},
\tag {3.33}
$$
by (3.23) and (3.26) we obtain
$$\displaystyle
\lim_{n\rightarrow\infty}\Vert\mbox{\boldmath $w$}_n(\,\cdot\,)-\mbox{\boldmath $w$}_x^*(\,\cdot\,;\mbox{\boldmath $a$})\Vert_{H^2(\Omega\setminus\overline{D})}=0.
$$
Now from this together with (3.31), (3.32) and (3.33), we obtain (3.28).

Next we describe the proof of (ii) and (iii).   Clearly we have
$$\displaystyle
I^*(x,\mbox{\boldmath $a$})\le
-\Vert 2\mu\,\text{Sym}\,\nabla(\mbox{\boldmath $J$}_{\Omega}(\,\cdot\,;x,\mbox{\boldmath $a$})\Vert_{L^2(D)}^2.
$$
By  (3.23) , (3.25) and the blowing up property already confirmed, that is, 
$$\displaystyle
\lim_{x\rightarrow a\in\partial D}\,\Vert\text{Sym}\,\nabla(\mbox{\boldmath $J$}(\,\cdot\,-x)\mbox{\boldmath $a$)}\Vert_{L^2(D)}^2=\infty,
$$
we conclude the validity of (3.29).

\noindent
Besides, for each positive numbers $\epsilon_1$ and $\epsilon_2$ it follows from (3.16), (3.25), (3.26)
and (3.27) that (3.30) is valid.

\noindent
$\Box$

The derivation of the next theorem is completely parallel to that of  (3.21)
except for a trick.

\proclaim{\noindent Theorem 3.7.}  Fix $\mbox{\boldmath $a$}\not=\mbox{\boldmath $0$}$.
Let $x\in\Omega\setminus\overline{D}$.
We have
$$\displaystyle
I^*(x,\mbox{\boldmath $a$})=\mbox{\boldmath $w$}_x^*(x;\mbox{\boldmath $a$})\cdot\mbox{\boldmath $a$}.
\tag {3.34}
$$
\endproclaim
{\it\noindent Proof.}
One can choose a needle $\sigma\in N_x$ with $\sigma\cap\overline{D}=\emptyset$.
Then there exists a needle sequence $\mbox{\boldmath $\xi$}=\{(\mbox{\boldmath $v$}_n,q_n)\}\in{\cal N}(\sigma,x,\mbox{\boldmath $a$})$.
Then, similar to Theorem 3.5, we have
$$\displaystyle
-\lim_{n\rightarrow\infty}
<(\Lambda_0-\Lambda_D)(\mbox{\boldmath $v$}_n
+\mbox{\boldmath $\epsilon$}_x(\,\cdot\,;\mbox{\boldmath $a$})\vert_{\partial\Omega},
(\mbox{\boldmath $J$}_{\Omega})_n(\,\cdot\,;x)\vert_{\partial\Omega}>
=\mbox{\boldmath $w$}_x^*(x;\mbox{\boldmath $a$})\cdot\mbox{\boldmath $a$},
\tag {3.35}
$$
$$\begin{array}{ll}
\displaystyle
(\mbox{\boldmath $J$}_{\Omega})_n(z;x)
=\mbox{\boldmath $J$}_{\Omega}(z;x)
-(\mbox{\boldmath $v$}_n(z)+\mbox{\boldmath $\epsilon$}_x(z;\mbox{\boldmath $a$})),
& z\in\Omega\setminus\{x\}.
\end{array}
$$
However, this right-hand takes the form
$$\begin{array}{l}
\displaystyle
\,\,\,\,\,\,
\mbox{\boldmath $J$}_{\Omega}(z;x)
-(\mbox{\boldmath $v$}_n(z)+\mbox{\boldmath $\epsilon$}_x(z;\mbox{\boldmath $a$}))
\\
\\
\displaystyle
=\mbox{\boldmath $J$}(z-x)\mbox{\boldmath $a$}+\mbox{\boldmath $\epsilon$}_x(z;\mbox{\boldmath $a$})
-(\mbox{\boldmath $v$}_n(z)+\mbox{\boldmath $\epsilon$}_x(z;\mbox{\boldmath $a$}))
\\
\\
\displaystyle
=\mbox{\boldmath $J$}(z-x)\mbox{\boldmath $a$}-\mbox{\boldmath $v$}_n(z)
\\
\\
\displaystyle
=\mbox{\boldmath $J$}_n(z;\mbox{\boldmath $\xi$}).
\end{array}
$$
Thus one gets
$$\displaystyle
(\mbox{\boldmath $J$}_{\Omega})_n(z;x)=\mbox{\boldmath $J$}_n(z;\mbox{\boldmath $\xi$}).
$$
And also we have, for $z\in\partial D$
$$\begin{array}{l}
\displaystyle
\,\,\,\,\,\,
\mbox{\boldmath $v$}_n(z)+\mbox{\boldmath $\epsilon$}_x(z;\mbox{\boldmath $a$})
\\
\\
\displaystyle
=\mbox{\boldmath $v$}_n(z)-\mbox{\boldmath $J$}(z-x)\mbox{\boldmath $a$}\\
\\
\displaystyle
=-\mbox{\boldmath $J$}_n(z;\mbox{\boldmath $\xi$}).
\end{array}
$$
From these we obtain
$$\begin{array}{l}
\displaystyle
\,\,\,\,\,\,
-<(\Lambda_0-\Lambda_D)(\mbox{\boldmath $v$}_n
+\mbox{\boldmath $\epsilon$}_x(\,\cdot\,;\mbox{\boldmath $a$}))\vert_{\partial\Omega},
(\mbox{\boldmath $J$}_{\Omega})_n(\,\cdot\,;x)\vert_{\partial\Omega}>
\\
\\
\displaystyle
=<(\Lambda_0-\Lambda_D)\mbox{\boldmath $J$}_n(\,\cdot\,;\mbox{\boldmath $\xi$})\vert_{\partial\Omega},
\mbox{\boldmath $J$}_n(\,\cdot\,;\mbox{\boldmath $\xi$})\vert_{\partial\Omega}>.
\end{array}
$$
This together with (3.28) and (3.35) yields the desired conclusion.

\noindent
$\Box$

$\quad$

{\bf\noindent Remark 3.2.}
Note that the proof presented above starts with the use of the {\it existence} of the needle sequence.
However, instead, using the representation formula of the reflected solution $\mbox{\boldmath $w$}_x^*(z;\mbox{\boldmath $a$})$ for $z\in\Omega\setminus\overline{D}$ the same as Theorem 3.3
one can obtain the formula (3.34).
We leave it for a reader as an exercise.

\subsubsection{Side B and unexpected result}

The following theorem deals with how the indicator sequence behaves when the tip of a needle enters or penetrates an obstacle.

\proclaim{\noindent Theorem 3.8.}  Fix $\mbox{\boldmath $a$}\not=\mbox{\boldmath $0$}$.
Let $x\in\Omega$ and $\sigma\in N_x$.  Assume that one of two cases (a) and (b) listed below are satisfied:

\noindent
(a) $x\in\overline{D}$;

\noindent
(b) $x\in\Omega\setminus\overline{D}$ and $\sigma\cap D\not=\emptyset$.

\noindent
Then, for any needle sequence $\mbox{\boldmath $\xi$}\in{\cal N}(x,\sigma,\mbox{\boldmath $a$})$ we have
$$\displaystyle
\lim_{n\rightarrow\infty}
<(\Lambda_0-\Lambda_D)\mbox{\boldmath $J$}_n(\,\cdot\,;\mbox{\boldmath $\xi$})\vert_{\partial\Omega},
\mbox{\boldmath $J$}_n(\,\cdot\,;\mbox{\boldmath $\xi$})\vert_{\partial\Omega}>
=-\infty.
$$

\endproclaim

As a direct corollary we obtain a characterization of  the set $\overline{D}$ 
in terms of the blowing up property of the sequence
$$\displaystyle
<(\Lambda_0-\Lambda_D)\mbox{\boldmath $J$}_n(\,\cdot\,;\mbox{\boldmath $\xi$})\vert_{\partial\Omega},
\mbox{\boldmath $J$}_n(\,\cdot\,;\mbox{\boldmath $\xi$})\vert_{\partial\Omega}>.
$$

$\quad$

\proclaim{\noindent Corollary 3.3. }  Fix $\mbox{\boldmath $a$}\not=\mbox{\boldmath $0$}$.
A point $x\in\Omega$ belongs to the set $\overline{D}$ if and only if
for any needle $\sigma\in N_x$ and needle sequence $\mbox{\boldmath $\xi$}\in
{\cal N}(x,\sigma,\mbox{\boldmath $a$})$
we have
$$\displaystyle
\lim_{n\rightarrow\infty}
<(\Lambda_0-\Lambda_D)\mbox{\boldmath $J$}_n(\,\cdot\,;\mbox{\boldmath $\xi$})\vert_{\partial\Omega},
\mbox{\boldmath $J$}_n(\,\cdot\,;\mbox{\boldmath $\xi$})\vert_{\partial\Omega}>
=-\infty.
$$
\endproclaim

And as an {\it unexpected byproduct}, we have a theorem that states the existence of the 
Side B for the singular sources method developed in subsection 3.2.

\proclaim{\noindent Theorem 3.9.}  Fix $\mbox{\boldmath $a$}\not=\mbox{\boldmath $0$}$.
Let $x\in\Omega$ and $\sigma\in N_x$.  Assume that one of two cases (a) and (b) listed below are satisfied:

\noindent
(a) $x\in\overline{D}$;

\noindent
(b) $x\in\Omega\setminus\overline{D}$ and $\sigma\cap D\not=\emptyset$.

\noindent
Then, for any needle sequence $\mbox{\boldmath $\xi$}=\{(\mbox{\boldmath $v$}_n,q_n)\}\in{\cal N}(x,\sigma,\mbox{\boldmath $a$})$ we have
$$\displaystyle
-\lim_{n\rightarrow\infty}
<(\Lambda_0-\Lambda_D)\mbox{\boldmath $v$}_n\vert_{\partial\Omega},
\mbox{\boldmath $J$}_n(\,\cdot\,;\mbox{\boldmath $\xi$})\vert_{\partial\Omega}>
=-\infty.
$$

\endproclaim

{\it\noindent Proof.}
Decompose
$$\begin{array}{ll}
\displaystyle
\mbox{\boldmath $J$}(z-x)\mbox{\boldmath $a$}
=\mbox{\boldmath $v$}_n(z)+\mbox{\boldmath $J$}_n(z;\mbox{\boldmath $\xi$}),
&
\displaystyle
z\in\partial\Omega.
\end{array}
$$
Then we have
$$\begin{array}{ll}
\displaystyle
<(\Lambda_0-\Lambda_D)\mbox{\boldmath $J$}(\,\cdot\,-x)\mbox{\boldmath $a$}\vert_{\partial\Omega},
\mbox{\boldmath $J$}(\,\cdot\,-x)\mbox{\boldmath $a$}\vert_{\partial\Omega}>
&
\displaystyle
=<(\Lambda_0-\Lambda_D)\mbox{\boldmath $v$}_n\vert_{\partial\Omega},
\mbox{\boldmath $v$}_n\vert_{\partial\Omega}>\\
\\
\displaystyle
&
\displaystyle
\,\,\,
+<(\Lambda_0-\Lambda_D)\mbox{\boldmath $J$}_n(\,\cdot\,;\mbox{\boldmath $\xi$})\vert_{\partial\Omega},
\mbox{\boldmath $J$}(\,\cdot\,;\mbox{\boldmath $\xi$})\vert_{\partial\Omega}>\\
\\
\displaystyle
&
\displaystyle
\,\,\,
+2<(\Lambda_0-\Lambda_D)\mbox{\boldmath $v$}_n\vert_{\partial\Omega},
\mbox{\boldmath $J$}_n(\,\cdot\,;\mbox{\boldmath $\xi$})\vert_{\partial\Omega}>.
\end{array}
$$
Thus one gets the representation
$$\begin{array}{l}
\,\,\,\,\,\,
\displaystyle
-<(\Lambda_0-\Lambda_D)\mbox{\boldmath $v$}_n\vert_{\partial\Omega},
\mbox{\boldmath $J$}_n(\,\cdot\,;\mbox{\boldmath $\xi$})\vert_{\partial\Omega}>
\\
\\
\displaystyle
=\frac{1}{2}
\left(<(\Lambda_0-\Lambda_D)\mbox{\boldmath $v$}_n\vert_{\partial\Omega},
\mbox{\boldmath $v$}_n\vert_{\partial\Omega}>
+<(\Lambda_0-\Lambda_D)\mbox{\boldmath $J$}_n(\,\cdot\,;\mbox{\boldmath $\xi$})\vert_{\partial\Omega},
\mbox{\boldmath $J$}(\,\cdot\,;\mbox{\boldmath $\xi$})\vert_{\partial\Omega}>
\right)\\
\\
\displaystyle
\,\,\,
-\frac{1}{2}<(\Lambda_0-\Lambda_D)\mbox{\boldmath $J$}(\,\cdot\,-x)\mbox{\boldmath $a$}\vert_{\partial\Omega},
\mbox{\boldmath $J$}(\,\cdot\,-x)\mbox{\boldmath $a$}\vert_{\partial\Omega}>.
\end{array}
\tag {3.36}
$$
Now the assertion of Theorem 3.9 follows from Theorems 3.2 and 3.8.

\noindent
$\Box$

As a direct corollory of Theorem 3.9 we immediately obtain

\proclaim{\noindent Corollary 3.4. }  Fix $\mbox{\boldmath $a$}\not=\mbox{\boldmath $0$}$.
A point $x\in\Omega$ belongs to the set $\overline{D}$ if and only if
for any needle $\sigma\in N_x$ and needle sequence $\mbox{\boldmath $\xi$}
=\{(\mbox{\boldmath $v$}_n,q_n)\}\in
{\cal N}(x,\sigma,\mbox{\boldmath $a$})$
we have
$$\displaystyle
-\lim_{n\rightarrow\infty}
<(\Lambda_0-\Lambda_D)\mbox{\boldmath $v$}_n\vert_{\partial\Omega},
\mbox{\boldmath $J$}_n(\,\cdot\,;\mbox{\boldmath $\xi$})\vert_{\partial\Omega}>
=-\infty.
$$
\endproclaim

\noindent
Needless to say, the facts corresponding to Theorem 3.8 and Corollary 3.4 also are valid
for the Laplace equation case \cite{IPS}.

Theorem 3.5, Theorem 3.9 and Corollary 3.4 tell us that the sequence
$$\displaystyle
\{-<(\Lambda_0-\Lambda_D)\mbox{\boldmath $v$}_n\vert_{\partial\Omega},
\mbox{\boldmath $J$}_n(\,\cdot\,;\mbox{\boldmath $\xi$})\vert_{\partial\Omega}>\}
$$
plays completely the same role as the indicator sequence for the Side B  of the probe method.
The Side B for the singular sources method exists!  This is the answer raised in \cite{IPS}.
Besides, the expression (3.36) means that its behaviour governed by the {\it mean-value} of
the indicator sequences for the Side B of the original and integrated probe methods.

\subsubsection{Relationship of all the indicator functions}

From the proof of Theorem 3.8 together with (3.5), (3.21) and (3.28) we obtain,
for all $x\in\Omega\setminus\overline{D}$
$$\displaystyle
\mbox{\boldmath $w$}_x(x)\mbox{\boldmath $a$}\cdot\mbox{\boldmath $a$}
=\frac{1}{2}(I(x,\mbox{\boldmath $a$})+I^*(x,\mbox{\boldmath $a$})
-<(\Lambda_0-\Lambda_D)\mbox{\boldmath $J$}(\,\cdot\,-x)\mbox{\boldmath $a$}\vert_{\partial\Omega},
\mbox{\boldmath $J$}(\,\cdot\,-x)\mbox{\boldmath $a$}\vert_{\partial\Omega}>).
\tag {3.37}
$$
Recalling (1.10), (3.4) and Definition 3.12, we see that this gives the {\it purely energy integral decomposition}
of the indicator function for the singular sources method formulated in Subsection 3.2.
Besides, it tells us that
the negative part of $\mbox{\boldmath $w$}_x(x)\mbox{\boldmath $a$}\cdot\mbox{\boldmath $a$}$
is given by
$$\displaystyle
\frac{1}{2}(I(x,\mbox{\boldmath $a$})+I^*(x,\mbox{\boldmath $a$}))
$$
and positive one given by
$$\displaystyle
-\frac{1}{2}
<(\Lambda_0-\Lambda_D)\mbox{\boldmath $J$}(\,\cdot\,-x)\mbox{\boldmath $a$}\vert_{\partial\Omega},
\mbox{\boldmath $J$}(\,\cdot\,-x)\mbox{\boldmath $a$}\vert_{\partial\Omega}>.
$$
Besides, (3.37) means that the indicator function of the singular sources method can be completely
expressed by using the indicator functions of  the original and integrated probe methods.

Finally we present two equations which describes the relationship between the all of the indicator functions.

\proclaim{\noindent Theorem 3.10.}  Let $x\in\Omega\setminus\overline{D}$.
It holds that
$$\begin{array}{ll}
\displaystyle
I^*(x,\mbox{\boldmath $a$})-I(x,\mbox{\boldmath $a$})
&
\displaystyle
=
2(I^1(x,\mbox{\boldmath $a$})-\mbox{\boldmath $w$}_x^1(x)\mbox{\boldmath $a$}\cdot\mbox{\boldmath $a$})
\\
\\
\displaystyle
&
\displaystyle
\,\,\,
+<(\Lambda_0-\Lambda_D)(\mbox{\boldmath $J$}(\cdot\,-x)\mbox{\boldmath $a$}\vert_{\partial\Omega},\mbox{\boldmath $J$}(\cdot\,-x)\mbox{\boldmath $a$}\vert_{\partial\Omega}>
\end{array}
\tag {3.38}
$$
and
$$\begin{array}{ll}
\displaystyle
\mbox{\boldmath $w$}^*_x(x;\mbox{\boldmath $a$})\cdot\mbox{\boldmath $a$}-\mbox{\boldmath $w$}_x(x)\mbox{\boldmath $a$}\cdot\mbox{\boldmath $a$}
&
\displaystyle
=
(I^1(x,\mbox{\boldmath $a$})-\mbox{\boldmath $w$}_x^1(x)\mbox{\boldmath $a$}\cdot\mbox{\boldmath $a$})
\\
\\
\displaystyle
&
\displaystyle
\,\,\,+
<(\Lambda_0-\Lambda_D)(\mbox{\boldmath $J$}(\cdot\,-x)\mbox{\boldmath $a$}\vert_{\partial\Omega},\mbox{\boldmath $J$}(\cdot\,-x)\mbox{\boldmath $a$}\vert_{\partial\Omega}>.
\end{array}
\tag {3.39}
$$

\endproclaim

\noindent
The proof of Theorem 3.10 is given in Subsection 4.4.
Note that the first term of  the right-hand sides on  (3.38) and (3.39), that is, 
$I^1(x,\mbox{\boldmath $a$})-\mbox{\boldmath $w$}_x^1(x)\mbox{\boldmath $a$}\cdot\mbox{\boldmath $a$}$
can be calculated from
$\Lambda_0-\Lambda_D$.  See (4.23).

\section{Proof of Theorems}

\subsection{Proof of Theorem 2.1}

We employ a modification of the argument for proofs of Lemmas 2.1 and 2.2 in \cite{INew}.

\noindent
(i)  Proof of (a). 
This is a consequence of the contradiction argument combined with the blowing up property of 
$\text{Sym}\nabla\mbox{\boldmath $J$}(y-x)\mbox{\boldmath $a$}$ described in Proposition 2.
The argument is parallel to the original argument in Lemma 2.1 of \cite{INew} and so we omit its
description.

\noindent
(ii) Proof of (b).   This part needs an essentially different argument from Lemma 2.2 in \cite{INew}.
Roughly speaking, the reason is, instead of the usual energy of $\nabla\mbox{\boldmath $v$}_n$
one must treat the strain energy $\text{Smy}\,\nabla\mbox{\boldmath $v$}_n$ which is a part of the former
energy.

Let $(\mbox{\boldmath $v$}, q)\in H^2(\Omega)\times H^1(\Omega)$ be an arbitrary pair of solutions of the Stokes system (1.8).
Since $\Delta q=0$ in $\Omega$, by the interior regularity,  $q$ is smooth in $\Omega$ and thus so is $\Delta\mbox{\boldmath $v$}$.  Thus $\mbox{\boldmath $v$}$ is also smooth in $\Omega$. 

{\bf\noindent Step 1.}
Let $B'$ be an open ball centered at $z$ such that $\overline{B'}\subset B\cap\Omega$.
Note that $B'$ can be arbitrary small if necessary.

{\bf\noindent Step 2.}
Choose a smaller open ball $B''$ centered at $z$ such that $\overline{B''}\subset B'$.

{\bf\noindent Step 3.}
Choose a $C^2$-domain $U$ in such a way that $\Sigma\equiv \partial U\cap\partial B''$ has a positive surface measure on $\partial B''$, $\text{dist}\,(\partial U\setminus\Sigma,\sigma)>0$, $x\in U$ and $\overline{U}\subset\Omega$.  

{\bf\noindent Step 4.}  Let $B'''$ be an open ball centered at $x$ such that $\overline{B'''}\subset U$.

\proclaim{\noindent Lemma 4.1.}
We have
$$\displaystyle
\Vert\mbox{\boldmath $v$}\Vert_{L^2(U)}
\le C
\Vert\mbox{\boldmath $v$}\Vert_{L^2(\partial U)}
\tag {4.1}
$$
and
$$\displaystyle
\Vert q\Vert_{L^2(U)}\le C\Vert q\Vert_{L^2(\partial U)}.
\tag {4.2}
$$
\endproclaim

{\it\noindent Proof.}
By the well-posedness for the Stokes system in $U$,  there exists a pair $(\mbox{\boldmath $z$}, r)\in H^2(U)\times H^1(U)$
such that
$$\left\{
\begin{array}{ll}
\displaystyle
\mu\Delta\mbox{\boldmath $z$}-\nabla r=\mbox{\boldmath $v$}, 
&
\displaystyle
y\in U,
\\
\\
\displaystyle
\nabla\cdot\mbox{\boldmath $z$}=0,
&
\displaystyle
y\in U\\
\\
\displaystyle
\mbox{\boldmath $z$}=\mbox{\boldmath $0$},
&
\displaystyle
y\in\partial U
\end{array}
\right.
$$
and 
$$\displaystyle
\Vert\mbox{\boldmath $z$}\Vert_{H^2(U)}+\Vert r-r_U\Vert_{H^1(U)}\le C\Vert\mbox{\boldmath $v$}\Vert_{L^2(U)},
\tag {4.3}
$$
where $C$ is a positive number independent of $\mbox{\boldmath $v$}$.
One may assume $r_U=0$ by a simple trick.

Integration by parts together with the governing equations of $(\mbox{\boldmath $z$}, r)$ and 
$(\mbox{\boldmath $v$},q)$ yields
$$\begin{array}{ll}
\displaystyle
\int_U\vert\mbox{\boldmath $v$}\vert^2\,dy
&
\displaystyle
=\int_{\partial U}(2\mu\,(\text{Sym}\,\nabla\mbox{\boldmath $z$})\mbox{\boldmath $\nu$}-r\mbox{\boldmath $\nu$})\cdot
\mbox{\boldmath $v$}\,dS(y).
\end{array}
$$
Applying  the trace theorem  and estimate (4.3) to this right-hand side, we obtain
$$\begin{array}
{ll}
\displaystyle
\Vert\mbox{\boldmath $v$}\Vert_{L^2(U)}^2
&
\displaystyle
\le (2\mu\Vert\text{Sym}\,\nabla\mbox{\boldmath $z$}\Vert_{L^2(\partial U)}+\Vert r\Vert_{L^2(\partial U)}\,)
\Vert\mbox{\boldmath $v$}\Vert_{L^2(\partial U)}
\\
\\
\displaystyle
&
\displaystyle
\le C(\Vert\mbox{\boldmath $z$}\Vert_{H^2(U)}+\Vert r\Vert_{H^1(U)}\,)\Vert\mbox{\boldmath $v$}\Vert_{L^2(\partial U)}
\\
\\
\displaystyle
&
\displaystyle
\le C\Vert\mbox{\boldmath $v$}\Vert_{L^2(U)}\Vert\mbox{\boldmath $v$}\Vert_{L^2(\partial U)}.
\end{array}
$$
Therefore we obtain (4.1).
Since the $q$ satisfies $\Delta q=0$ in $U$, a similar and simpler argument yields also (4.2) (see also, for example, Lemma A.2 in \cite{INew}).

\noindent
$\Box$

Next we prepare the following estimate.

\proclaim{\noindent Lemma 4.2.}
We have
$$\displaystyle
\int_{B'''}\vert\text{Sym}\,\nabla\mbox{\boldmath $v$}\vert^2\,dy
\le C\left(\int_{U}\vert\mbox{\boldmath $v$}\vert^2\,dy
+\int_{U}\vert q\vert^2\,dy\right).
\tag {4.4}
$$

\endproclaim

{\it\noindent Proof.}
Let $\varphi\in C_0^{\infty}(U)$ and multiply the equation
$\mu\Delta\mbox{\boldmath $v$}-\nabla q=\mbox{\boldmath $0$}$ in $U$ by $\varphi^2\mbox{\boldmath $v$}$.
Then, we have
$$\begin{array}{ll}
\displaystyle
0
&
\displaystyle
=\int_U\left(2\mu\,\text{Sym}\nabla\mbox{\boldmath $v$}-qI_3)\right)\cdot\nabla(\varphi^2\mbox{\boldmath $v$})\,dy
\\
\\
\displaystyle
&
\displaystyle
=\int_U2\mu\,\text{Sym}\nabla\mbox{\boldmath $v$}\cdot\,\text{Sym}\,\nabla(\varphi^2\mbox{\boldmath $v$})\,dy
-\int_Uq\nabla\cdot(\varphi^2\mbox{\boldmath $v$})\,dy
\\
\\
\displaystyle
&
\displaystyle
=
\int_U4\mu\,\varphi\text{Sym}\,\nabla\mbox{\boldmath $v$}\cdot\,\text{Sym}\,(\mbox{\boldmath $v$}\otimes\nabla\varphi)\,dy
+\int_U2\mu\varphi^2\vert\text{Sym}\,(\nabla\mbox{\boldmath $v$})\vert^2\,dy
\\
\\
\displaystyle
&
\displaystyle
\,\,\,
-\int_U2q\,\varphi\nabla\varphi\cdot\mbox{\boldmath $v$}\,dy.
\end{array}
\tag {4.5}
$$
From (4.5) we have
$$\begin{array}{ll}
\displaystyle
\int_U2\mu\varphi^2\vert\text{Sym}\,\nabla\mbox{\boldmath $v$}\vert^2\,dy
&
\displaystyle
=-\int_U4\mu\,\varphi\,\text{Sym}\,\nabla\mbox{\boldmath $v$}\cdot\,\text{Sym}\,(\mbox{\boldmath $v$}\otimes\nabla\varphi)\,dy
+\int_U2q\,\varphi\nabla\varphi\cdot\mbox{\boldmath $v$}\,dy.
\end{array}
$$
Rewrite this as
$$
\displaystyle
\int_U2\mu\vert \varphi\text{Sym}\,\nabla\mbox{\boldmath $v$}+\text{Sym}\,(\mbox{\boldmath $v$}\otimes\nabla\varphi)\vert^2\,dy
=\int_U2\mu\vert\text{Sym}\,(\mbox{\boldmath $v$}\otimes\nabla\varphi)\vert^2\,dy
+\int_U 2q\varphi\nabla\varphi\cdot\mbox{\boldmath $v$}\,dy.
\tag {4.6}
$$
From (4.6) we have
$$\displaystyle
\frac{1}{2}\int_U2\mu\varphi^2\vert\text{Sym}\,\nabla\mbox{\boldmath $v$}\vert^2\,dy
\le 
2\int_U2\mu\vert\text{Sym}\,(\mbox{\boldmath $v$}\otimes\nabla\varphi)\vert^2\,dy
+\int_U2\vert q\vert\vert\varphi\vert\nabla\varphi\vert\vert\mbox{\boldmath $v$}\vert\,dy,
$$
i.e.,
$$\displaystyle
\int_U2\mu\varphi^2\vert\text{Sym}\,\nabla\mbox{\boldmath $v$}\vert^2\,dy
\le 
4\int_U2\mu\vert\text{Sym}\,(\mbox{\boldmath $v$}\otimes\nabla\varphi)\vert^2\,dy
+2\int_U2\vert q\vert\vert\varphi\vert\nabla\varphi\vert\vert\mbox{\boldmath $v$}\vert\,dy.
\tag {4.7}
$$
Now choose $\varphi$ in such a way that $\varphi=1$ on $\overline{B'''}$.
Then, from (4.7) we obtain
$$\begin{array}{ll}\displaystyle
\int_{B'''}\vert\text{Sym}\,\nabla\mbox{\boldmath $v$}\vert^2\,dy
&
\displaystyle
\le C\left(\int_U\vert\mbox{\boldmath $v$}\vert^2\,dy
+\int_{U}\vert q\vert\vert\mbox{\boldmath $v$}\vert\,dy.
\right).
\end{array}
$$
This yields (4.4).

\noindent
$\Box$

Then, by Lemma 4.1 and (4.4), we obtain
$$\displaystyle
\int_{B'''}\vert\text{Sym}\,\nabla\mbox{\boldmath $v$}\vert^2\,dy
\le C\left(\int_{\partial U}\vert\mbox{\boldmath $v$}\vert^2\,dS(y)
+\int_{\partial U}\vert q\vert^2\,dS(y)\right).
\tag {4.8}
$$

Let us continue the proof of (b).
Divide the integrals on the  right-hand side of  (4.8) as
$$\displaystyle
\int_{\partial U}\vert\mbox{\boldmath $v$}\vert^2\,dS(y)
+\int_{\partial U}\vert q\vert^2\,dS(y)
=\int_{\Sigma}\vert\mbox{\boldmath $v$}\vert^2\,dS(y)
+\int_{\Sigma}\vert q\vert^2\,dS(y)
+{\cal R}(\mbox{\boldmath $v$}, q),
$$
where
$$\displaystyle
{\cal R}(\mbox{\boldmath $v$}, q)
=\int_{\partial U\setminus\Sigma}\vert\mbox{\boldmath $v$}\vert^2\,dS(y)
+\int_{\partial U\setminus\Sigma}\vert q\vert^2\,dS(y).
$$
We have
$$\begin{array}{ll}
\displaystyle
\int_{\Sigma}\vert\mbox{\boldmath $v$}\vert^2\,dS(y)
+\int_{\Sigma}\vert q\vert^2\,dS(y)
&
\displaystyle
\le \int_{\partial B''}\vert\mbox{\boldmath $v$}\vert^2\,dS(y)
+\int_{\partial B''}\vert q\vert^2\,dS(y)
\\
\\
\displaystyle
&
\displaystyle
\le C
(\Vert\mbox{\boldmath $v$}\Vert_{H^1(B'')}^2
+\Vert q\Vert_{H^1(B'')}^2).
\end{array}
$$
Here applying Korn's second inequality to the first term on this right-hand side, we obtain
$$\begin{array}{ll}
\displaystyle
\int_{\Sigma}\vert\mbox{\boldmath $v$}\vert^2\,dS(y)
+\int_{\Sigma}\vert q\vert^2\,dS(y)
&
\displaystyle
\le
C(\Vert \mbox{\boldmath $v$}\Vert_{L^2(B'')}^2+\Vert\text{Sym}\,\nabla\mbox{\boldmath $v$}\Vert_{L^2(B'')}^2
+\Vert q\Vert_{H^1(B'')}^2).
\end{array}
$$
This together with (4.8) yields
$$\displaystyle
\Vert\text{Sym}\,\nabla\mbox{\boldmath $v$}\Vert_{L^2(B''')}^2
\le C(\Vert \mbox{\boldmath $v$}\Vert_{L^2(B'')}^2+\Vert\text{Sym}\,\nabla\mbox{\boldmath $v$}\Vert_{L^2(B'')}^2
+\Vert q\Vert_{H^1(B'')}^2+{\cal R}(\mbox{\boldmath $v$},q)).
\tag {4.9}
$$
Here we consider how to estimate $\Vert q\Vert_{H^1(B'')}^2$ in terms of $\mbox{\boldmath $v$}$ over $B'$
modulo a convergent term.
We employ a variant of the Poincar\'e lemma (see, e.g., \cite{SS} and \cite{Z}).

\proclaim{\noindent Lemma 4.3.}
Let $E\subset B''$ be an arbitrary measurable subset of $B''$ with $\vert E\vert>0$.
We have
$$\displaystyle
\Vert q-q_E\Vert_{L^2(B'')}\le C_E\Vert\nabla q\Vert_{L^2(B'')}.
$$

\endproclaim

Applying this lemma to the term $\Vert q\Vert_{H^1(B'')}^2$,
we have
$$\displaystyle
 \Vert q\Vert_{H^1(B'')}^2
 \le C_E(\Vert\nabla q\Vert_{L^2(B'')}^2+\vert q_E\vert^2\vert B''\vert).
 $$
Here the equation $\mu\Delta\mbox{\boldmath $v$}=\nabla q$  and $\nabla\cdot\mbox{\boldmath $v$}=0$
yields
$$\displaystyle
\Vert\nabla q\Vert_{L^2(B'')}=\mu\Vert\Delta\mbox{\boldmath $v$}\Vert_{L^2(B'')}
=\mu\Vert\nabla\times(\nabla\times\mbox{\boldmath $v$})\Vert_{L^2(B'')}.
$$
Since $\Delta(\nabla\times\mbox{\boldmath $v$})=\mbox{\boldmath $0$}$, a similar argument as Lemma 4.1,
we have
$$\begin{array}{ll}
\displaystyle
\Vert\nabla\times(\nabla\times\mbox{\boldmath $v$})\Vert_{L^2(B'')}
&
\displaystyle
\le C\Vert \nabla\times\mbox{\boldmath $v$}\Vert_{L^2(B')}
\\
\\
\displaystyle
&
\displaystyle
\le C'\Vert\mbox{\boldmath $v$}\Vert_{H^1(B')}
\end{array}
$$
Thus Korn's second inequality \cite{DuL} yields
$$\displaystyle
\Vert\nabla\times(\nabla\times\mbox{\boldmath $v$})\Vert_{L^2(B'')}^2
\le C(\Vert\mbox{\boldmath $v$}\Vert_{L^2(B')}^2+\Vert\text{Sym}\nabla\mbox{\boldmath $v$}\Vert_{L^2(B')}^2\,).
$$
Therefore we obtain
$$\displaystyle
\Vert q\Vert_{H^1(B'')}^2
\le
C_E(\Vert\mbox{\boldmath $v$}\Vert_{L^2(B')}^2+\Vert\text{Sym}\,\nabla\mbox{\boldmath $v$}\Vert_{L^2(B')}^2+\vert q_E\vert^2\vert B''\vert).
$$
This together with (4.9) yields
$$
\displaystyle
\Vert\text{Sym}\,\nabla\mbox{\boldmath $v$}\Vert_{L^2(B''')}^2
\le C(\Vert \mbox{\boldmath $v$}\Vert_{L^2(B')}^2+\Vert\text{Sym}\,\nabla\mbox{\boldmath $v$}\Vert_{L^2(B')}^2
+
{\cal R}(\mbox{\boldmath $v$},q)+\vert q_E\vert^2\vert B''\vert).
\tag {4.10}
$$

Here we give an upper estimate for $\Vert \mbox{\boldmath $v$}\Vert_{L^2(B')}^2$ of the right-hand side on (4.11).
We make use of  a variant of Korn's Poincar\'e-type inequality (which is a corollary of \cite{N}, Theorem 2.2, see also \cite{BP}, Corollary 1.2)
$$\displaystyle
\Vert\mbox{\boldmath $v$}-\Pi_{B'}\mbox{\boldmath $v$}\Vert_{L^2(B')}
\le C_{B'}\Vert\text{Sym}\,\nabla\mbox{\boldmath $v$}\Vert_{L^2(B')},
\tag {4.11}
$$
where $\Pi_{B'}:L^2(B')\rightarrow {\cal R}$ denotes the orthogonal projection of $L^2(B')$ onto 
the set of all rigid displacement fields\footnote{Note that we have another expression of ${\cal R}$:
$$\displaystyle
{\cal R}=\{\,d(y)\mbox{\boldmath $h$}\,\vert\, \mbox{\boldmath $h$}\in\Bbb R^6\,\}.
$$
}
$$\displaystyle
{\cal R}=\{\mbox{\boldmath $b$}+\mbox{\boldmath $c$}\times x\,\vert\, \mbox{\boldmath $b$}, \mbox{\boldmath $c$}\in\Bbb R^3\,\}\subset L^2(B').
$$
For our purpose, we need the explicit expression of $\Pi_{B'}$.  By \cite{BP}, 
we have
$$\displaystyle
\Pi_{B'}(\mbox{\boldmath $v$})(y)
=d(y)<\mbox{\boldmath $v$}>_{B'},
$$
where
$$\displaystyle
d(y)=\left(
\begin{array}{cc}
\displaystyle
I_3 
&
\displaystyle
\tilde{d}(y)
\end{array}
\right),
$$
$$\displaystyle
\tilde{d}(x)
=
\left(
\begin{array}{ccc}
\displaystyle
0
&
\displaystyle
-y_3
&
\displaystyle
y_2
\\
\\
\displaystyle
y_3
&
\displaystyle
0
&
\displaystyle
-y_1
\\
\\
\displaystyle
-y_2
&
\displaystyle
y_1
&
\displaystyle
0
\end{array}
\right),
$$
$$\displaystyle
d_{B'}=\int_{B'}d(y)^Td(y)\,dy
$$
and
$$\displaystyle
<\mbox{\boldmath $v$}>_{B'}
=d_{B'}^{-1}\int_{B'}d(y)^T\mbox{\boldmath $v$}(y)\,dy.
$$

Let $E\subset B'$ be an arbitrary nonempty open subset of $B'$.
Then
$$\displaystyle
\Vert\mbox{\boldmath $v$}-\Pi_{E}\mbox{\boldmath $v$}\Vert_{L^2(B')}
\le \Vert\mbox{\boldmath $v$}-\Pi_{B'}\mbox{\boldmath $v$}\Vert_{L^2(B')}
+\Vert\Pi_E\mbox{\boldmath $v$}-\Pi_{B'}\mbox{\boldmath $v$}\Vert_{L^2(B')}.
\tag {4.12}
$$
Here we show that
$$
\displaystyle
\Vert\Pi_{E}\mbox{\boldmath $v$}-\Pi_{B'}\mbox{\boldmath $v$}\Vert_{L^2(B')}
\le
C_{B',E}\Vert\mbox{\boldmath $v$}-\Pi_{B'}\mbox{\boldmath $v$}\Vert_{L^2(B')}.
\tag {4.13}
$$
Here we employ an argument done in the proof of Proposition 2.3 in \cite{INew} which is taken from \cite{SS}.
Write
$$\begin{array}{ll}
\displaystyle
\Pi_E\mbox{\boldmath $v$}
&
\displaystyle
=d<\mbox{\boldmath $v$}>_E\\
\\
\displaystyle
&
\displaystyle
=dd_E^{-1}\int_E d^T\mbox{\boldmath $v$}\,dy
\\
\\
\displaystyle
&
\displaystyle
=dd_E^{-1}\int_E d^T(\mbox{\boldmath $v$}-d<\mbox{\boldmath $v$}>_{B'})\,dy
+dd_E^{-1}\int_Ed^Td<\mbox{\boldmath $v$}>_{B'}\,dy
\\
\\
\displaystyle
&
\displaystyle
=dd_E^{-1}\int_E d^T(\mbox{\boldmath $v$}-d<\mbox{\boldmath $v$}>_{B'})\,dy
+d<\mbox{\boldmath $v$}>_{B'}
\\
\\
\displaystyle
&
\displaystyle
=dd_E^{-1}\int_E d^T(\mbox{\boldmath $v$}-\Pi_{B'}\mbox{\boldmath $v$})\,dy
+\Pi_{B'}\mbox{\boldmath $v$}.
\end{array}
$$
Thus one gets
$$
\begin{array}{ll}
\displaystyle
\Vert\Pi_{E}\mbox{\boldmath $v$}-\Pi_{B'}\mbox{\boldmath $v$}\Vert_{L^2(B')}
&
\displaystyle
=\Vert dd_E^{-1}\int_E d^T(\mbox{\boldmath $v$}-\Pi_{B'}\mbox{\boldmath $v$})\,dy\Vert_{L^2(B')}
\\
\\
\displaystyle
&
\displaystyle
=\Vert dd_E^{-1}\Vert_{L^2(B')}
\left\vert
\int_E d^T(\mbox{\boldmath $v$}-\Pi_{B'}\mbox{\boldmath $v$}\,)dy\right\vert
\\
\\
\displaystyle
&
\displaystyle
\le
\Vert dd_E^{-1}\Vert_{L^2(B')}
\Vert d^T\Vert_{L^2(E)}
\Vert\mbox{\boldmath $v$}-\Pi_{B'}\mbox{\boldmath $v$}\Vert_{L^2(B')}.
\end{array}
$$
Therefore (4.13) is valid for $C_{B',E}=\Vert dd_E^{-1}\Vert_{L^2(B')}
\Vert d^T\Vert_{L^2(E)}$.
Thus from (4.12) and (4.13) one gets
$$\displaystyle
\Vert\mbox{\boldmath $v$}-\Pi_{E}\mbox{\boldmath $v$}\Vert_{L^2(B')}
\le (1+C_{B',E})\Vert\mbox{\boldmath $v$}-\Pi_{B'}\mbox{\boldmath $v$}\Vert_{L^2(B')}.
\tag {4.14}
$$
This together with (4.11) yields
$$\begin{array}{ll}
\displaystyle
\Vert \mbox{\boldmath $v$}\Vert_{L^2(B')}
&
\displaystyle
\le\Vert \mbox{\boldmath $v$}-\Pi_{E}\mbox{\boldmath $v$}\Vert_{L^2(B')}+
\Vert\Pi_E\mbox{\boldmath $v$}\Vert_{L^2(B')}
\\
\\
\displaystyle
&
\displaystyle
\le C_{B'}(1+C_{B', E})\Vert\text{Sym}\,\nabla\mbox{\boldmath $v$}\Vert_{L^2(B')}
+\Vert\Pi_E\mbox{\boldmath $v$}\Vert_{L^2(B')}.
\end{array}
$$
Thus this together with (4.10) yields
$$\begin{array}{ll}
\displaystyle
\Vert\text{Sym}\nabla\mbox{\boldmath $v$}\Vert_{L^2(B''')}^2
&
\displaystyle
\le
C(B,E')
\left(\Vert\text{Sym}\nabla\mbox{\boldmath $v$}\Vert_{L^2(B')}^2
\right.
\\
\\
\displaystyle
&
\displaystyle
\,\,\,
\left.
+\Vert\Pi_{E}\mbox{\boldmath $v$}\Vert_{L^2(B')}^2
+{\cal R}(\mbox{\boldmath $v$},q)+\vert q_E\vert^2\vert B''\vert
\,\right).
\end{array}
\tag {4.15}
$$

Once we have (4.15), the proof of (b) is as follows.
Choose $E$ such that $\overline{E}\cap\sigma=\emptyset$.
This choice is always possible since $B'\setminus\sigma\not=\emptyset$ and is open.   
Substituting $\mbox{\boldmath $v$}=\mbox{\boldmath $v$}_n$ into (4.15), we have
$$\begin{array}{ll}
\displaystyle
\Vert\text{Sym}\nabla\mbox{\boldmath $v$}_n\Vert_{L^2(B')}^2
&
\displaystyle
\ge
C(B',E)^{-1}
\Vert\text{Sym}\nabla\mbox{\boldmath $v$}_n\Vert_{L^2(B''')}^2
\\
\\
\displaystyle
&
\displaystyle
\,\,\,
-\Vert\Pi_{E}\mbox{\boldmath $v$}_n\Vert_{L^2(B')}^2-{\cal R}(\mbox{\boldmath $v$}_n,q_n)
-\vert (q_n)_E\vert^2\vert B''\vert.
\end{array}
$$
Using the convergence property of the needle sequence on $\partial U\setminus\Sigma$, one sees the boundedness of  the term ${\cal R}(\mbox{\boldmath $v$}_n,q_n)$.
And also $\Vert\Pi_E\mbox{\boldmath $v$}_n\Vert_{L^2(B')}$ is bounded by the choice of $E$
and the explicit formula of $\Pi_E$, and so is $\vert (q_n)_E\vert^2$.
Thus from the blowing up property of $\Vert\text{Sym}\,\nabla\mbox{\boldmath $v$}_n\Vert_{L^2(B''')}^2$
which is a consequence of  (a), one can conclude the blowing up
of  $\Vert\text{Sym}\nabla\mbox{\boldmath $v$}_n\Vert_{L^2(B')}^2$
 on $B'$ and thus $B\cap\Omega$, too.
This completes the proof of (b) in Theorem 2.1.

\subsection{Proof of Theorem 3.3}

We start with the standard integral representation formulae  of $\mbox{\boldmath $w$}_x$ and $\mbox{\boldmath $w$}_x^1$, which are the consequences of the divergence theorem.

\proclaim{\noindent Proposition 4.1}
Let $\mbox{\boldmath $a$}$ and $\mbox{\boldmath $b$}$ be arbitrary real constant vectors.
We have, for all $(x,y)\in(\Omega\setminus\overline{D})\times(\Omega\setminus\overline{D})$
$$\begin{array}{l}
\displaystyle
\,\,\,\,\,\,
\mbox{\boldmath $w$}_x(y)\mbox{\boldmath $a$}\cdot\mbox{\boldmath $b$}
\\
\\
\displaystyle
=\int_{\partial\Omega}\sigma(\mbox{\boldmath $w$}_x(z)\mbox{\boldmath $a$},q_x(z))\mbox{\boldmath $\nu$}\cdot
\mbox{\boldmath $J$}(z-y)\mbox{\boldmath $b$}\,dz
\\
\\
\displaystyle
\,\,\,
-
\int_{\partial D}
\mbox{\boldmath $J$}(z-x)\mbox{\boldmath $a$}\cdot
\sigma(\mbox{\boldmath $J$}(z-y)\mbox{\boldmath $b$},\mbox{\boldmath $p$}(z-y)\cdot\mbox{\boldmath $b$}))\mbox{\boldmath $\nu$}
\,dz\\
\\
\displaystyle
\,\,\,
-\int_{\partial D}
\sigma(\mbox{\boldmath $w$}_x(z)\mbox{\boldmath $a$},q_x(z))\mbox{\boldmath $\nu$}\cdot
\mbox{\boldmath $J$}(z-y)\mbox{\boldmath $b$}\,dz
\end{array}
\tag {4.16}
$$
and
$$\begin{array}{l}
\displaystyle
\,\,\,\,\,\,
\mbox{\boldmath $w$}_x^1(y)\mbox{\boldmath $a$}\cdot\mbox{\boldmath $b$}
\\
\\
\displaystyle
=\int_{\partial\Omega}\sigma(\mbox{\boldmath $w$}_x^1(z)\mbox{\boldmath $a$},q_x^1(z))\mbox{\boldmath $\nu$}\cdot
\mbox{\boldmath $J$}(z-y)\mbox{\boldmath $b$}\,dz\\
\\
\displaystyle
\,\,\,-\int_{\partial\Omega}\mbox{\boldmath $J$}(z-x)\mbox{\boldmath $a$}\cdot
\sigma(\mbox{\boldmath $J$}(z-y)\mbox{\boldmath $b$},\mbox{\boldmath $p$}(z-y)\cdot\mbox{\boldmath $b$})\mbox{\boldmath $\nu$}\,dS(z)
\\
\\
\displaystyle
\,\,\,
-\int_{\partial D}\sigma(\mbox{\boldmath $w$}_x^1(z)\mbox{\boldmath $a$},q_x^1(z))\mbox{\boldmath $\nu$}\cdot
\mbox{\boldmath $J$}(z-y)\mbox{\boldmath $b$}\,dz.
\end{array}
\tag {4.17}
$$

\endproclaim

{\it\noindent Proof.}
First we represent $\mbox{\boldmath $w$}_x(y)\mbox{\boldmath $a$}$, $y\in\Omega\setminus\overline{D}$ in terms of
the Cauchy data of $\mbox{\boldmath $w$}_x\mbox{\boldmath $a$}$ on $\partial\Omega$ and $\partial D$.
We have
$$\begin{array}{l}
\displaystyle
\,\,\,\,\,\,
\mbox{\boldmath $w$}_x(y)\mbox{\boldmath $a$}\cdot\mbox{\boldmath $b$}
\\
\\
\displaystyle
=\int_{\Omega\setminus\overline{D}}\mbox{\boldmath $w$}_x(z)\mbox{\boldmath $a$}\cdot\mbox{\boldmath $b$}\,\delta(z-y)\,dz
\\
\\
\displaystyle
=-\int_{\Omega\setminus\overline{D}}\mbox{\boldmath $w$}_x(z)\mbox{\boldmath $a$}\cdot
\left\{\mu\Delta(\mbox{\boldmath $J$}(z-y)\mbox{\boldmath $b$})-\nabla(\mbox{\boldmath $p$}(z-y)\cdot\mbox{\boldmath $b$})\right\}\,dz
\\
\\
\displaystyle
=-\int_{\Omega\setminus\overline{D}}\mbox{\boldmath $w$}_x(z)\mbox{\boldmath $a$}\cdot
\text{div}\,\sigma(\mbox{\boldmath $J$}(z-y)\mbox{\boldmath $b$},\mbox{\boldmath $p$}(z-y)\cdot\mbox{\boldmath $b$})\,dz
\\
\\
\displaystyle
=-\int_{\partial D}\mbox{\boldmath $J$}(z-x)\mbox{\boldmath $a$}\cdot
\sigma(\mbox{\boldmath $J$}(z-y)\mbox{\boldmath $b$},
\mbox{\boldmath $p$}(z-y)\cdot\mbox{\boldmath $b$})\mbox{\boldmath $\nu$}\,dS(z)
\\
\\
\displaystyle
\,\,\,
+\int_{\Omega\setminus\overline{D}}\nabla\mbox{\boldmath $w$}_x(z)\mbox{\boldmath $a$}
\cdot\sigma(\mbox{\boldmath $J$}(z-y)\mbox{\boldmath $b$},
\mbox{\boldmath $p$}(z-y)\cdot\mbox{\boldmath $b$}))\,dz.
\end{array}
\tag {4.18}
$$
Since $\nabla(\mbox{\boldmath $w$}_x\mbox{\boldmath $a$})\cdot I_3=\nabla\cdot(\mbox{\boldmath $w$}_x\mbox{\boldmath $a$})=\nabla\cdot
 (\mbox{\boldmath $J$}(z-y)\mbox{\boldmath $b$})=0$ in
$\Omega\setminus\overline{D}$, we have
$$\begin{array}{l}
\displaystyle
\,\,\,\,\,\,
\int_{\Omega\setminus\overline{D}}\nabla(\mbox{\boldmath $w$}_x(z)\mbox{\boldmath $a$})
\cdot\sigma(\mbox{\boldmath $J$}(z-y)\mbox{\boldmath $b$},
\mbox{\boldmath $p$}(z-y)\cdot\mbox{\boldmath $b$}))\,dz
\\
\\
\displaystyle
=\int_{\Omega\setminus\overline{D}}
2\mu\nabla(\mbox{\boldmath $w$}_x(z)\mbox{\boldmath $a$})\cdot\text{Sym}\,\nabla(\mbox{\boldmath $J$}(z-y)\mbox{\boldmath $b$})
\,dz
\\
\\
\displaystyle
=\int_{\Omega\setminus\overline{D}}
\sigma(\mbox{\boldmath $w$}_x(z)\mbox{\boldmath $a$},q_x(z))\cdot\nabla(\mbox{\boldmath $J$}(z-y)\mbox{\boldmath $b$})
\,dz
\\
\\
\displaystyle
=\int_{\partial\Omega}\sigma(\mbox{\boldmath $w$}_x(z)\mbox{\boldmath $a$},q_x(z))\mbox{\boldmath $\nu$}\cdot
\mbox{\boldmath $J$}(z-y)\mbox{\boldmath $b$}\,dz
-\int_{\partial D}\sigma(\mbox{\boldmath $w$}_x(z)\mbox{\boldmath $a$},q_x(z))\mbox{\boldmath $\nu$}\cdot
\mbox{\boldmath $J$}(z-y)\mbox{\boldmath $b$}\,dz.
\end{array}
$$
Substituting this into (4.18), we obtain (4.16).

Similarly we have
$$\begin{array}{l}
\displaystyle
\,\,\,\,\,\,
\mbox{\boldmath $w$}_x^1(y)\mbox{\boldmath $a$}\cdot\mbox{\boldmath $b$}
\\
\\
\displaystyle
=\int_{\Omega\setminus\overline{D}}\mbox{\boldmath $w$}_x^1(z)\cdot\mbox{\boldmath $b$}\,\delta(z-y)\,dz
\\
\\
\displaystyle
=-\int_{\Omega\setminus\overline{D}}\mbox{\boldmath $w$}_x^1(z)\mbox{\boldmath $a$}\cdot
\left\{\mu\Delta(\mbox{\boldmath $J$}(z-y)\mbox{\boldmath $b$})-\nabla p(z-y;\mbox{\boldmath $b$})\right\}\,dz
\\
\\
\displaystyle
=-\int_{\Omega\setminus\overline{D}}\mbox{\boldmath $w$}_x^1(z)\mbox{\boldmath $a$}\cdot
\text{div}\,\sigma(\mbox{\boldmath $J$}(z-y)\mbox{\boldmath $b$},p(z-y;\mbox{\boldmath $b$})\,dz
\\
\\
\displaystyle
=-\int_{\partial\Omega}\mbox{\boldmath $J$}(z-x)\mbox{\boldmath $a$}\cdot
\sigma(\mbox{\boldmath $J$}(z-y)\mbox{\boldmath $b$},p(z-y;\mbox{\boldmath $b$})\mbox{\boldmath $\nu$}\,dS(z)
\\
\\
\displaystyle
\,\,\,
+\int_{\Omega\setminus\overline{D}}\nabla(\mbox{\boldmath $w$}_x^1(z)\mbox{\boldmath $a$})
\cdot\sigma(\mbox{\boldmath $J$}(z-y)\mbox{\boldmath $b$},p(z-y;\mbox{\boldmath $b$}))\,dz.
\end{array}
\tag {4.19}
$$
Since $\nabla(\mbox{\boldmath $w$}_x^1\mbox{\boldmath $a$})\cdot I_3=\nabla\cdot(\mbox{\boldmath $w$}_x^1\mbox{\boldmath $a$})=\nabla\cdot
 (\mbox{\boldmath $J$}(z-y)\mbox{\boldmath $b$})=0$ in
$\Omega\setminus\overline{D}$, we have
$$\begin{array}{l}
\displaystyle
\,\,\,\,\,\,
\int_{\Omega\setminus\overline{D}}\nabla(\mbox{\boldmath $w$}_x^1(z)\mbox{\boldmath $a$})
\cdot\sigma(\mbox{\boldmath $J$}(z-y)\mbox{\boldmath $b$},p(z-y;\mbox{\boldmath $b$}))\,dz
\\
\\
\displaystyle
=\int_{\Omega\setminus\overline{D}}
2\mu\nabla(\mbox{\boldmath $w$}_x^1(z)\mbox{\boldmath $a$})\cdot\text{Sym}\,\nabla(\mbox{\boldmath $J$}(z-y)\mbox{\boldmath $b$})
\,dz
\\
\\
\displaystyle
=\int_{\Omega\setminus\overline{D}}
\sigma(\mbox{\boldmath $w$}_x^1(z)\mbox{\boldmath $a$},q_x^1(z;\mbox{\boldmath $a$}))\cdot\nabla(\mbox{\boldmath $J$}(z-y)\mbox{\boldmath $b$})
\,dz
\\
\\
\displaystyle
=\int_{\partial\Omega}\sigma(\mbox{\boldmath $w$}_x^1(z)\mbox{\boldmath $a$},q_x^1(z;\mbox{\boldmath $a$}))\mbox{\boldmath $\nu$}\cdot
\mbox{\boldmath $J$}(z-y)\mbox{\boldmath $b$}\,dz
-\int_{\partial D}\sigma(\mbox{\boldmath $w$}_x^1(z)\mbox{\boldmath $a$},
q_x^1(z;\mbox{\boldmath $a$}))\mbox{\boldmath $\nu$}\cdot
\mbox{\boldmath $J$}(z-y)\mbox{\boldmath $b$}\,dz.
\end{array}
$$
Substituting this into (4.19), we obtain (4.17).

\noindent
$\Box$

We further rewrite (4.16) and (4.17) as follows.

\proclaim{\noindent Proposition 4.2.}
We have, for all $(x,y)\in (\Omega\setminus\overline{D})^2$
$$\begin{array}{ll}
\displaystyle
\,\,\,\,\,\,
\mbox{\boldmath $w$}_x(y)\mbox{\boldmath $a$}\cdot\mbox{\boldmath $b$}
&
\displaystyle
=\int_{\partial\Omega}\sigma(\mbox{\boldmath $w$}_x(z)\mbox{\boldmath $a$},q_x(z;\mbox{\boldmath $a$}))\mbox{\boldmath $\nu$}\cdot
\mbox{\boldmath $J$}(z-y)\mbox{\boldmath $b$}\,dz
\\
\\
\displaystyle
&
\displaystyle
\,\,\,
-\int_D2\mu\,\text{Sym}\,\nabla(\mbox{\boldmath $J$}(z-x)\mbox{\boldmath $a$})
\cdot
\text{Sym}\,\nabla(\mbox{\boldmath $J$}(z-y)\mbox{\boldmath $b$})
\,dz
\\
\\
\displaystyle
&
\displaystyle
\,\,\,
-\int_{\Omega\setminus\overline{D}}2\mu\,\text{Sym}\nabla(\mbox{\boldmath $w$}_x(z)\mbox{\boldmath $a$})
\cdot\text{Sym}\nabla(\mbox{\boldmath $w$}_y(z)\mbox{\boldmath $b$})\,dz
\end{array}
\tag {4.20}
$$
and
$$\begin{array}{ll}
\displaystyle
\mbox{\boldmath $w$}_x^1(y)\mbox{\boldmath $a$}\cdot\mbox{\boldmath $b$}
&
\displaystyle
=\int_{\Omega\setminus\overline{D}}2\mu\,\text{Sym}\nabla(\mbox{\boldmath $w$}_x^1(z)\mbox{\boldmath $a$})
\cdot\text{Sym}\nabla(\mbox{\boldmath $w$}_y^1(z)\mbox{\boldmath $b$})\,dz
\\
\\
\displaystyle
&
\displaystyle
\,\,\,
+
\int_{\Bbb R^3\setminus\overline{\Omega}}
2\mu\,\text{Sym}\,\nabla(\mbox{\boldmath $J$}(z-x)\mbox{\boldmath $a$})
\cdot\text{Sym}\,\nabla(\mbox{\boldmath $J$}(z-y)\mbox{\boldmath $b$})\,dz
\\
\\
\displaystyle
&
\displaystyle
\,\,\,
-\int_{\partial\Omega}\mbox{\boldmath $J$}(z-x)\mbox{\boldmath $a$}
\cdot\sigma(\mbox{\boldmath $w$}_y(z)\mbox{\boldmath $b$},q_y(z;\mbox{\boldmath $b$}))
\mbox{\boldmath $\nu$}\,dS(z).
\end{array}
\tag {4.21}
$$

\endproclaim

{\it\noindent Proof.}
We have
$$\begin{array}{l}
\displaystyle
\,\,\,\,\,\,
-\int_{\partial D}
\sigma(\mbox{\boldmath $w$}_x(z)\mbox{\boldmath $a$},q_x(z;\mbox{\boldmath $a$}))\mbox{\boldmath $\nu$}\cdot
\mbox{\boldmath $J$}(z-y)\mbox{\boldmath $b$}\,dz
\\
\\
\displaystyle
=\int_{\partial D}
\sigma(\mbox{\boldmath $w$}_x(z)\mbox{\boldmath $a$},q_x(z;\mbox{\boldmath $a$}))\mbox{\boldmath $\nu$}\cdot
\mbox{\boldmath $w$}_y(z)\mbox{\boldmath $b$}\,dz
\\
\\
\displaystyle
=-
\left(
\int_{\partial\Omega}\sigma(\mbox{\boldmath $w$}_x(z)\mbox{\boldmath $a$},q_x(z;\mbox{\boldmath $a$}))\mbox{\boldmath $\nu$}\cdot
\mbox{\boldmath $w$}_y(z)\mbox{\boldmath $b$}\,dz
-\int_{\partial D}\sigma(\mbox{\boldmath $w$}_x(z)\mbox{\boldmath $a$},q_x(z;\mbox{\boldmath $a$}))\mbox{\boldmath $\nu$}\cdot
\mbox{\boldmath $w$}_y(z)\mbox{\boldmath $b$}\,dz
\right)
\\
\\
\displaystyle
=-\int_{\Omega\setminus\overline{D}}
\text{div}(\sigma(\mbox{\boldmath $w$}_x(z)\mbox{\boldmath $a$},q_x(z;\mbox{\boldmath $a$}))\cdot
\mbox{\boldmath $w$}_y(z)\mbox{\boldmath $b$}\,dz\\
\\
\displaystyle
\,\,\,
-\int_{\Omega\setminus\overline{D}}\sigma(\mbox{\boldmath $w$}_x(z)\mbox{\boldmath $a$},q_x(z;\mbox{\boldmath $a$}))
\cdot\nabla(\mbox{\boldmath $w$}_y(z)\mbox{\boldmath $b$})\,dz
\\
\\
\displaystyle
=-\int_{\Omega\setminus\overline{D}}2\mu\,\text{Sym}\nabla(\mbox{\boldmath $w$}_x(z)\mbox{\boldmath $a$})
\cdot\text{Sym}\nabla(\mbox{\boldmath $w$}_y(z)\mbox{\boldmath $b$})\,dz
\end{array}
$$
and 
$$\begin{array}{l}
\displaystyle
\,\,\,\,\,\,
\int_{\partial D}
\mbox{\boldmath $J$}(z-x)\mbox{\boldmath $a$}\cdot
\sigma(\mbox{\boldmath $J$}(z-y)\mbox{\boldmath $b$},
\mbox{\boldmath $p$}(z-y)\cdot\mbox{\boldmath $b$}))\mbox{\boldmath $\nu$}
\,dz
\\
\\
\displaystyle
=\int_{D}\mbox{\boldmath $J$}(z-x)\mbox{\boldmath $a$}\cdot
\text{div}\,(\sigma(\mbox{\boldmath $J$}(z-y)\mbox{\boldmath $b$},
\mbox{\boldmath $p$}(z-y)\cdot\mbox{\boldmath $b$}))\,dz
\\
\\
\displaystyle
\,\,\,
+\int_D2\mu\,\text{Sym}\,\nabla(\mbox{\boldmath $J$}(z-x)\mbox{\boldmath $a$})
\cdot
\text{Sym}\,\nabla(\mbox{\boldmath $J$}(z-y)\mbox{\boldmath $b$})
\,dz
\\
\\
\displaystyle
=\int_D2\mu\,\text{Sym}\,\nabla(\mbox{\boldmath $J$}(z-x)\mbox{\boldmath $a$})
\cdot
\text{Sym}\,\nabla(\mbox{\boldmath $J$}(z-y)\mbox{\boldmath $b$})
\,dz.
\end{array}
$$
Substituting these into (4.16) we obtain (4.20).

Next we have
$$\begin{array}{l}
\displaystyle
\,\,\,\,\,\,
\int_{\partial\Omega}\sigma(\mbox{\boldmath $w$}_x^1(z)\mbox{\boldmath $a$},q_x^1(z;\mbox{\boldmath $a$}))\mbox{\boldmath $\nu$}\cdot
\mbox{\boldmath $J$}(z-y)\mbox{\boldmath $b$}\,dz\\
\\
\displaystyle
=\int_{\partial\Omega}\sigma(\mbox{\boldmath $w$}_x^1(z)\mbox{\boldmath $a$},q_x^1(z;\mbox{\boldmath $a$}))\mbox{\boldmath $\nu$}
\cdot
\mbox{\boldmath $w$}_y^1(z)\mbox{\boldmath $b$}\,dz
\\
\\
\displaystyle
\,\,\,
-\int_{\partial D}\sigma(\mbox{\boldmath $w$}_x^1(z)\mbox{\boldmath $a$},q_x^1(z;\mbox{\boldmath $a$}))\mbox{\boldmath $\nu$}
\cdot
\mbox{\boldmath $w$}_y^1(z)\mbox{\boldmath $b$}\,dz
\\
\\
\displaystyle
=\int_{\Omega\setminus\overline{D}}\text{div}(\sigma(\mbox{\boldmath $w$}_x^1(z)\mbox{\boldmath $a$},q_x^1(z;\mbox{\boldmath $a$})))
\cdot
\mbox{\boldmath $w$}_y^1(z)\mbox{\boldmath $b$}\,dz
\\
\\
\displaystyle
\,\,\,
+\int_{\Omega\setminus\overline{D}}2\mu\,\text{Sym}\nabla(\mbox{\boldmath $w$}_x^1(z)\mbox{\boldmath $a$})
\cdot\text{Sym}\nabla(\mbox{\boldmath $w$}_y^1(z)\mbox{\boldmath $b$})\,dz
\\
\\
\displaystyle
=\int_{\Omega\setminus\overline{D}}2\mu\,\text{Sym}\nabla(\mbox{\boldmath $w$}_x^1(z)\mbox{\boldmath $a$})
\cdot\text{Sym}\nabla(\mbox{\boldmath $w$}_y^1(z)\mbox{\boldmath $b$})\,dz.
\end{array}
$$
Besides,
$$\begin{array}{l}
\displaystyle
\,\,\,\,\,\,
-\int_{\partial\Omega}\mbox{\boldmath $J$}(z-x)\mbox{\boldmath $a$}\cdot
\sigma(\mbox{\boldmath $J$}(z-y)\mbox{\boldmath $b$},
\mbox{\boldmath $p$}(z-y)\cdot\mbox{\boldmath $b$})\mbox{\boldmath $\nu$}\,dS(z)
\\
\\
\displaystyle
=\int_{\Bbb R^3\setminus\overline{\Omega}}
\mbox{\boldmath $J$}(z-x)\mbox{\boldmath $a$}\cdot
\text{div}\,(\sigma(\mbox{\boldmath $J$}(z-y)\mbox{\boldmath $b$},
\mbox{\boldmath $p$}(z-y)\cdot\mbox{\boldmath $b$}))\,dS(z)\\
\\
\displaystyle
\,\,\,
+\int_{\Bbb R^3\setminus\overline{\Omega}}
2\mu\,\text{Sym}\,\nabla(\mbox{\boldmath $J$}(z-x)\mbox{\boldmath $a$})
\cdot\text{Sym}\,\nabla(\mbox{\boldmath $J$}(z-y)\mbox{\boldmath $b$})
\,dz
\\
\\
\displaystyle
=\int_{\Bbb R^3\setminus\overline{\Omega}}
2\mu\,\text{Sym}\,\nabla(\mbox{\boldmath $J$}(z-x)\mbox{\boldmath $a$})
\cdot\text{Sym}\,\nabla(\mbox{\boldmath $J$}(z-y)\mbox{\boldmath $b$})\,dz
\end{array}
$$
and
$$\begin{array}{l}
\displaystyle
-\int_{\partial D}\sigma(\mbox{\boldmath $w$}_x^1(z)\mbox{\boldmath $a$},q_x^1(z;\mbox{\boldmath $a$}))\mbox{\boldmath $\nu$}\cdot
\mbox{\boldmath $J$}(z-y)\mbox{\boldmath $b$}\,dz
\\
\\
\displaystyle
=\int_{\partial D}\sigma(\mbox{\boldmath $w$}_x^1(z)\mbox{\boldmath $a$},q_x^1(z;\mbox{\boldmath $a$}))\mbox{\boldmath $\nu$}\cdot
\mbox{\boldmath $w$}_y(z)\mbox{\boldmath $b$}\,dz
\\
\\
\displaystyle
=-
\left(\int_{\partial\Omega}\sigma(\mbox{\boldmath $w$}_x^1(z)\mbox{\boldmath $a$},q_x^1(z;\mbox{\boldmath $a$}))\mbox{\boldmath $\nu$}\cdot
\mbox{\boldmath $w$}_y(z)\mbox{\boldmath $b$}\,dz
-\int_{\partial D}\sigma(\mbox{\boldmath $w$}_x^1(z)\mbox{\boldmath $a$},q_x^1(z;\mbox{\boldmath $a$}))\mbox{\boldmath $\nu$}\cdot
\mbox{\boldmath $w$}_y(z)\mbox{\boldmath $b$}\,dz
\right)
\\
\\
\displaystyle
=-\int_{\Omega\setminus\overline{D}}\text{div}\,(\sigma(\mbox{\boldmath $w$}_x^1(z)\mbox{\boldmath $a$},q_x^1(z;\mbox{\boldmath $a$})))
\cdot\mbox{\boldmath $w$}_y(z)\mbox{\boldmath $b$}\,dz
-\int_{\Omega\setminus\overline{D}}\sigma(\mbox{\boldmath $w$}_x^1(z)\mbox{\boldmath $a$},q_x^1(z;\mbox{\boldmath $a$}))
\cdot\nabla(\mbox{\boldmath $w$}_y(z)\mbox{\boldmath $b$})\,dz
\\
\\
\displaystyle
=-\int_{\Omega\setminus\overline{D}}\sigma(\mbox{\boldmath $w$}_x^1(z)\mbox{\boldmath $a$},q_x^1(z;\mbox{\boldmath $a$}))
\cdot\nabla(\mbox{\boldmath $w$}_y(z)\mbox{\boldmath $b$})\,dz
\\
\\
\displaystyle
=-\int_{\partial\Omega}\mbox{\boldmath $w$}_x^1(z)\mbox{\boldmath $a$}
\cdot\sigma(\mbox{\boldmath $w$}_y(z)\mbox{\boldmath $b$},q_y(z;\mbox{\boldmath $b$})
\mbox{\boldmath $\nu$}\,dS(z)
\\
\\
\displaystyle
=-\int_{\partial\Omega}\mbox{\boldmath $J$}(z-x)\mbox{\boldmath $a$}
\cdot\sigma(\mbox{\boldmath $w$}_y(z)\mbox{\boldmath $b$},q_y(z;\mbox{\boldmath $b$}))
\mbox{\boldmath $\nu$}\,dS(z).
\end{array}
$$
Substituting these into (4.17), we obtain (4.21).

\noindent
$\Box$

Now the expression (3.13) of $\mbox{\boldmath $W$}_x(y)\mbox{\boldmath $a$}\cdot\mbox{\boldmath $b$}$
in Theorem 3.3 is a direct consequence of (4.20) and (4.21).

\subsection{Proof of Theorem 3.5}

Write
$$\begin{array}{l}
\displaystyle
\,\,\,\,\,\,
-<(\Lambda_0-\Lambda_D)\mbox{\boldmath $v$}_n\vert_{\partial\Omega},
\mbox{\boldmath $J$}_n(\,\cdot\,;\mbox{\boldmath $\xi$})\vert_{\partial\Omega}>
\\
\\
\displaystyle
=-<(\Lambda_0-\Lambda_D)\mbox{\boldmath $v$}_n\vert_{\partial\Omega},
\mbox{\boldmath $J$}(\,\cdot\,-x)\mbox{\boldmath $a$})\vert_{\partial\Omega}>
+<(\Lambda_0-\Lambda_D)\mbox{\boldmath $v$}_n\vert_{\partial\Omega},
\mbox{\boldmath $v$}\vert_{\partial\Omega}>.
\end{array}
\tag {4.22}
$$
From (3.3) of Proposition 3.1 and (4.21) of  Proposition  4.2, we have
$$\displaystyle
I^1(x,\mbox{\boldmath $a$})-\mbox{\boldmath $w$}_x^1(x)\mbox{\boldmath $a$}\cdot\mbox{\boldmath $a$}
=-\lim_{n\rightarrow\infty}\,<(\Lambda_0-\Lambda_D)\mbox{\boldmath $v$}_n\vert_{\partial\Omega}, \mbox{\boldmath $J$}(\,\cdot\,-x)\mbox{\boldmath $a$}\vert_{\partial\Omega}>.
\tag {4.23}
$$
Now from (4.22), (4.23) and (3.5) of Theorem 3.1 yields
$$\displaystyle
-\lim_{n\rightarrow\infty}<(\Lambda_0-\Lambda_D)\mbox{\boldmath $v$}_n\vert_{\partial\Omega},
\mbox{\boldmath $J$}_n(\,\cdot\,;\mbox{\boldmath $\xi$})\vert_{\partial\Omega}>
=I^1(x,\mbox{\boldmath $a$})-\mbox{\boldmath $w$}_x^1(x)\mbox{\boldmath $a$}\cdot\mbox{\boldmath $a$}
+I(x,\mbox{\boldmath $a$}).
$$
Then the decomposition (3.11) and (3.14) tells us that this right-hand side coincides with $\mbox{\boldmath $w$}_x(x)\mbox{\boldmath $a$}\cdot\mbox{\boldmath $a$}$.

\subsection{Proof of Theorem 3.10.}

Here we present a proof of making use of the needle sequence, which is not elementary in the sense
that it is a consequence of the Runge approximation property.

There exists a needle $\sigma\in N_x$ satisfying $\sigma\cap\overline{D}=\emptyset$.
Let $\mbox{\boldmath $\xi$}=\{(\mbox{\boldmath $v$}_n,q_n)\}\in{\cal N}(x,\sigma,\mbox{\boldmath $a$})$.
Using the first equation on (3.20), one can rewrite
$$\begin{array}{l}
\displaystyle
\,\,\,\,\,\,
<(\Lambda_0-\Lambda_D)\mbox{\boldmath $J$}_n(\,\cdot\,;\mbox{\boldmath $\xi$})\vert_{\partial\Omega},
\mbox{\boldmath $J$}_n(\,\cdot\,;\mbox{\boldmath $\xi$})\vert_{\partial\Omega}>
\\
\\
\displaystyle
=<(\Lambda_0-\Lambda_D)(\mbox{\boldmath $J$}(\,\cdot\,-x)\mbox{\boldmath $a$}-\mbox{\boldmath $v$}_n)\vert_{\partial\Omega},
(\mbox{\boldmath $J$}(\,\cdot\,-x)\mbox{\boldmath $a$}-\mbox{\boldmath $v$}_n)\vert_{\partial\Omega}>
\\
\\
\displaystyle
=<(\Lambda_0-\Lambda_D)\mbox{\boldmath $v$}_n\vert_{\partial\Omega},
\mbox{\boldmath $v$}_n\vert_{\partial\Omega}>\\
\\
\displaystyle
\,\,\,
-<(\Lambda_0-\Lambda_D)\mbox{\boldmath $J$}(\,\cdot\,-x)\mbox{\boldmath $a$}\vert_{\partial\Omega},
\mbox{\boldmath $v$}_n\vert_{\partial\Omega}>
-<(\Lambda_0-\Lambda_D)\mbox{\boldmath $v$}_n\vert_{\partial\Omega},
\mbox{\boldmath $J$}(\,\cdot\,-x)\mbox{\boldmath $a$}\vert_{\partial\Omega}>
\\
\\
\displaystyle
\,\,\,
+
<(\Lambda_0-\Lambda_D)(\mbox{\boldmath $J$}(\cdot\,-x)\mbox{\boldmath $a$}\vert_{\partial\Omega},\mbox{\boldmath $J$}(\cdot\,-x)\mbox{\boldmath $a$}\vert_{\partial\Omega}>
\\
\\
\displaystyle
=<(\Lambda_0-\Lambda_D)\mbox{\boldmath $v$}_n\vert_{\partial\Omega},
\mbox{\boldmath $v$}_n\vert_{\partial\Omega}>
\\
\\
\displaystyle
\,\,\,
-2<(\Lambda_0-\Lambda_D)\mbox{\boldmath $v$}_n\vert_{\partial\Omega},
\mbox{\boldmath $J$}(\,\cdot\,-x)\mbox{\boldmath $a$}\vert_{\partial\Omega}>
\\
\\
\displaystyle
\,\,\,
+
<(\Lambda_0-\Lambda_D)(\mbox{\boldmath $J$}(\cdot\,-x)\mbox{\boldmath $a$}\vert_{\partial\Omega},\mbox{\boldmath $J$}(\cdot\,-x)\mbox{\boldmath $a$}\vert_{\partial\Omega}>.
\end{array}
$$
Then (4.23) together with (3.5) of  Theorem 3.1 and  (3.28) of Theorem 3.6 yields (3.38).

The equation (3.39) is a direct consequence of rewriting equation (3.38) via the equations
(3.11),  (3.14) of Theorem 3.3  and (3.34) of Theorem 3.7.

\end{document}